\documentclass{article}
\usepackage{graphicx} 
\usepackage[english]{babel}
\usepackage{amsfonts,amsthm,amsmath,amssymb,mathtools,nicefrac}
\usepackage{fullpage}
\usepackage{authblk}
\usepackage{enumitem}
\usepackage{todonotes}
\usepackage{hyperref}
\usepackage{bbm,stmaryrd}
\usepackage[most]{tcolorbox}
\usepackage{xr}
\newcommand{\By}[2]{\overset{\mbox{\tiny{#1}}}{#2}}
\newcommand{\ByRef}[2]{   \By{\eqref{#1}}{#2} }
\newcommand{\eqBy}[1]{    \By{#1}{=} }

\newcommand{\eqByRef}[1]{ \ByRef{#1}{=} }

\newcommand{\geByRef}[1]{ \ByRef{#1}{\ge} }
\newcommand{\JUSTIFY}[1]{\fbox{\tiny{#1}}\quad}

\newtheorem{lem}{Lemma}[section]		
\newtheorem{fact}[lem]{Fact}		
		
\newtheorem{defi}[lem]{Definition}
\newtheorem{rem}[lem]{Remark}
\newtheorem{thm}[lem]{Theorem}
\newtheorem{cor}[lem]{Corollary}		
		
\newtheorem{prop}[lem]{Proposition}		
		
\newtheorem{claim}{Claim}

\newcommand{\PSpec}{\mathrm{PSpec}} 
\newcommand{\Spec}{\mathrm{Spec}}

\newcommand{\N}{\mathbb{N}}
\newcommand{\Z}{\mathbb{Z}}
\newcommand{\R}{\mathbb{R}}
\newcommand{\C}{\mathbb{C}}
\newcommand{\diff}{\mathsf{d}}%differential symbol
\newcommand{\eps}{\varepsilon}
\newcommand{\ImaginaryUnit}{\mathbf{i}}
 
\newcommand{\Identity}{\mathrm{Id}} %Identity on a Banach space

\newcommand{\Krein}{Kre\u{\i}n}

\DeclareMathOperator{\esssup}{ess\, sup}
\DeclareMathOperator{\essinf}{ess\, inf}

\DeclareMathOperator{\degminIn}{mindeg^{in}}
\DeclareMathOperator{\degminOut}{mindeg^{out}}
\DeclareMathOperator{\degmaxIn}{maxdeg^{in}}
\DeclareMathOperator{\degmaxOut}{maxdeg^{out}}
\DeclareMathOperator{\degIn}{deg^{in}}
\DeclareMathOperator{\degOut}{deg^{out}}
\DeclareMathOperator{\ReachIn}{Reach^{in}}
\DeclareMathOperator{\ReachOut}{Reach^{out}}
\DeclareMathOperator{\ground}{ground}
\DeclareMathOperator{\supp}{supp}
\newcommand{\image}{\mathrm{image}}
\DeclareMathOperator{\cutnormdist}{d_{\square}}
\DeclareMathOperator{\cutdist}{\delta_{\square}}

\title{Digraphons: connectivity and spectral aspects}
\author{Jan Hladký, Petr Savický}
\date{}
\affil{Institute of Computer Science of the Czech Academy of Sciences\thanks{Research supported by Czech Science Foundation Project 21-21762X. With institutional support RVO:67985807.}}

\begin{document}
\maketitle
\begin{abstract}
The theory of graphons has proven to be a powerful tool in many areas of graph theory. In this paper, we introduce several foundational aspects of the theory of digraphons --- asymmetric two-variable functions that arise as limits of sequences of directed graphs (digraphs). Our results address their decomposition into strongly connected components, periodicity, spectral properties, and asymptotic behaviour of their large powers.
\end{abstract}

\section{Introduction}\label{sec:Intro}
The theory of limits of dense graph sequences, developed by Borgs, Chayes, Lovász, Sós, Szegedy, and Vesztergombi~\cite{LOVASZ2006933,BORGS20081801}, has proven to be a powerful framework in modern combinatorics. A typical application begins with a sequence of finite graphs and uses the Compactness Theorem of Lovász and Szegedy~\cite{MR2306658} to transfer a given problem to the analytic setting of graphons. Graphons, which are measurable functions representing limits of dense graphs, provide a versatile and robust language for addressing a wide range of problems. The richer the analytic tools available for handling graphons, the stronger and more flexible the theory becomes. Throughout the paper, given a probability space $(\Omega,\mu)$ (with an implicit sigma-algebra), a \emph{graphon} is a measurable function $W:\Omega^2\to [0,1]$ which is symmetric (that is, $W(x,y)=W(y,x)$).

The original theory of graph limits focuses on simple, undirected finite graphs. Since its inception, the scope of the theory has expanded significantly, encompassing limits of various combinatorial structures, including hypergraphs of fixed uniformity~\cite{ELEK20121731}, permutations~\cite{HOPPEN201393}, and Latin squares~\cite{Garbe2023Limits}, among others. All these limit theories retain core features of dense graph limits, such as compactness of the limit space and homomorphism densities, but they vary in how closely they adhere to the broader toolkit of graph limit theory. Even within the realm of graphs, several natural extensions have been studied. When the assumption of simplicity is relaxed to allow uniformly bounded edge multiplicities, the theory remains essentially unchanged. However, permitting unbounded edge multiplicities introduces new challenges and leads to more intricate limit objects, as explored in~\cite{KUNSZENTIKOVACS2022109284}.

In this paper, we focus on directed graphs (digraphs), where each pair of vertices $\{u,v\}$ may exhibit one of four configurations: no directed edge, a single directed edge $(u,v)$, a single directed edge $(v,u)$, or a bidirectional edge pair $(u,v)$ and $(v,u)$. Many core elements of the dense graph limit theory extend to digraphs \emph{mutatis mutandis}. The limit objects in this setting, called \emph{digraphons}, are naturally defined: just as the adjacency matrix of a digraph need not be symmetric, a digraphon is simply a graphon without the symmetry condition, that is, any measurable function $\Gamma:\Omega^2\to[0,1]$. The analog of the Compactness Theorem holds as well, with the main technical difference being the use of a variant of the Regularity Lemma adapted for digraphs. Many other key notions need only minimal adaptations. This includes the notion of homomorphism densities.

The focus of this work is to develop tools for digraphons that differ significantly from the graphon setting.

We arrived at this topic while studying a variant of the classical random 2-SAT problem. In this variant, the inclusion probabilities of individual clauses over variables $\{v_1, \ldots, v_n\}$ are governed by a graphon $W$, which serves as a parameter of the model. The central result of the corresponding paper~\cite{HladkySavicky:Inhomo2SAT} identifies a critical threshold $\rho^*(W)$, such that the behavior of the random formula transitions from being asymptotically almost surely satisfiable when $\rho^*(W) < 1$ to unsatisfiable when $\rho^*(W) > 1$. Although the original model uses a symmetric graphon to reflect the symmetry between clauses $(l_i \vee l_j)$ and $(l_j \vee l_i)$, the analysis naturally evolves into an asymmetric framework of digraphons. In fact, $\rho^*(W)$ turns out to be a spectral parameter associated with a digraphon derived from $W$. This led us to realize that many fundamental properties of digraphons had not been previously studied. Thus, several of the general results developed in this work were in fact originally motivated by the technical demands of~\cite{HladkySavicky:Inhomo2SAT}, but clearly stand as contributions of broad interest. Previous work on digraphons includes~\cite{cai2016priors}, which develops the theory of exchangeability in the directed setting and the corresponding Bayesian nonparametric statistical framework. Limits of directed graphs were also successfully employed to address extremal problems. Several papers in this direction, such as~\cite{MR4467141,MR3638333} do not employ digraphons but rather the more abstract framework of Razborov's flag algebras. We note that among the three digraphon notions discussed in this paper --- connectivity, periodicity, and spectrum --- the first two can be formulated within the framework of flag algebras, albeit not in a particularly explicit manner, whereas the third cannot. Digraphons were used for extremal analysis in~\cite{GRZESIK2023117,zhao2020impartial,MR4498587}. As a matter of fact, these three papers deal with limits of tournaments, which lead to particular digraphons, called `tournamentons' in~\cite{GRZESIK2023117} and `tourneyons' in \cite{zhao2020impartial}. For us, paper~\cite{GRZESIK2023117} is particularly relevant as it features prominently spectral techniques, which are also the main focus of our paper. A certain auxiliary step in~\cite{GRZESIK2023117} was handled in a somewhat artificial way. Very recently, Grebík, Kr\'al', Liu, Pikhurko, and Slipantschuk~\cite{grebik2025convergencespectradigraphlimits} proved that the eigenvalues of digraphons are continuous in the cut distance. This statement allows replacing the artificial step in~\cite{GRZESIK2023117} with a natural and compact argument.

Last, we mention~\cite{THORNBLAD201896}, which deals with the decomposition of tournamentons into components. Our results about the decomposition of digraphons into strong components generalize this result.

\section{Our results}

\subsection{Connectivity}
In Definition \ref{def:component} below, we introduce a collection of related connectivity concepts for digraphons. While for graphons, connectivity was introduced by Janson~\cite{ConnectednessGraphons} and extensively studied since, for digraphons these notions are new. To motivate them with the finite setting, recall that a nonempty set of vertices $C$ in digraph is \emph{strongly connected} if for every partition $A\sqcup B=C$ into two nonempty sets, there is at least one directed edge going from $A$ to $B$. Maximal strongly connected sets are called \emph{strong components}. The last notion in Definition~\ref{def:component} is of `fragmented sets'. This concept does not make any sense in finite digraphs. That is, one way to decompose a digraph into strong components is to start with an initial one-cell partition of the entire vertex set. If there is at any moment a cell violating the above condition on strong connectedness, we subdivide that cell accordingly. In a finite graph, this process must eventually terminate (since single vertices cannot be subdivided), resulting in the unique decomposition into connected components. In a digraphon, the regions where this process continues to split sets of positive measure into ever smaller ones constitute its fragmented sets.

\begin{figure}\centering
	\includegraphics[scale=0.7]{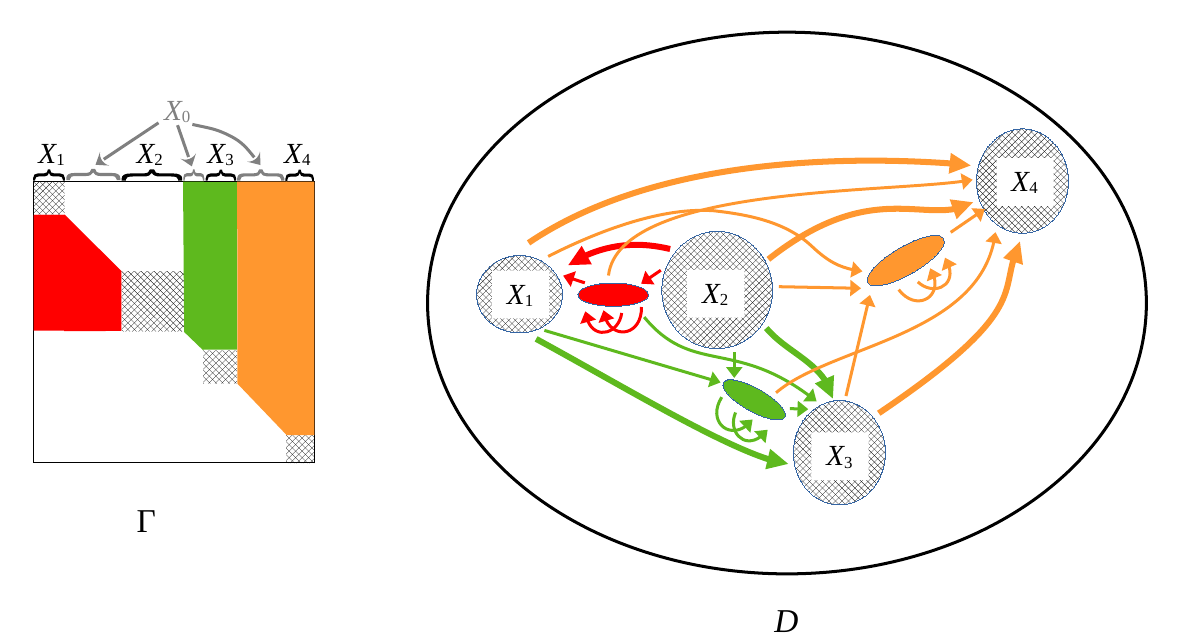}
	\caption{An example of a digraphon $\Gamma$ with its decomposition into strong components $X_1,\ldots,X_4$ and the remaining fragmented set $X_0$. Arbitrary strongly connected parts can be placed in the hatched parts. The white part is zero while the coloured parts are arbitrary (zero or positive). The digraphon $\Gamma$ should be thought of as a limit representation of a large digraph $D$. The three coloured sets are (subgraphs of) three transitive tournaments.}
	\label{fig:StrongComponents}
\end{figure}

\begin{defi}[strongly connected set, strong component, fragmented set]\label{def:component}
Suppose that $\Gamma$ is a digraphon on a probability space $(\Omega,\mu)$, and a set $X\subseteq \Omega$ is of positive measure.
\begin{enumerate}[label=(\roman*)]
    \item\label{en:StrongConn} We say that \emph{$X$ is strongly connected in $\Gamma$} if for every partition $A\sqcup B=X$ with $\mu(A),\mu(B)>0$ we have that $\int_{A\times B}\Gamma>0$.
    \item\label{en:StrongComp} We say that \emph{$X$ is a strong component in $\Gamma$} if $X$ is strongly connected and for every $Y\subseteq \Omega$ with $\mu(X\cap Y)>0$ and $\mu(Y\setminus X)>0$ we have that $Y$ is not strongly connected.
    \item\label{en:Fragmented} We say that \emph{$X$ is fragmented in $\Gamma$} if every subset $Y\subseteq X$ of positive measure is not strongly connected.
\end{enumerate}
\end{defi}
Recall that each digraph can be decomposed in a unique way into maximal strong components. Theorem~\ref{thm:decompositionIntoComponents} below is a digraphon counterpart. See Figure~\ref{fig:StrongComponents} for an illustration.

\begin{thm}\label{thm:decompositionIntoComponents}
Suppose that $\Gamma$ is a digraphon on $\Omega$. Then there exists a finite or a countable set $I$ not containing~$0$ and a decomposition $\Omega=X_0\sqcup \bigsqcup_{i\in I} X_i$ so that $X_0$ is either an empty set or is fragmented in $\Gamma$ and each $X_i$ is a strong component. 

Further, this decomposition is unique in the sense that if partitions $\{X_i\}_{i\in I\cup\{0\}}$ and $\{X'_i\}_{i\in I'\cup\{0\}}$ are two decompositions of $\Gamma$ into strong components as above, then there exists a bijection $\pi:I\to I'$ such that $X'_{\pi(i)}$ equals $X_i$ modulo a nullset for each $i\in I$, and $X_0$ equals $X'_0$ modulo a nullset.
\end{thm}
We prove Theorem~\ref{thm:decompositionIntoComponents} in Section~\ref{sec:ProofDecomposition}. In Section~\ref{ssec:AdditionalDigraphDecomposition} we establish further properties of decompositions into strong components. For example, these properties imply, that there is a natural digraphon counterpart of the notion of `condensation digraph' known from decompositions of digraphs into strong components.

Theorem~\ref{thm:decompositionIntoComponents} generalizes one of the main results of~\cite{THORNBLAD201896}. More specifically, Theorem~6.11 of~\cite{THORNBLAD201896} proves Theorem~\ref{thm:decompositionIntoComponents} for `tournamentons', that is, for digraphons $\Gamma$ satisfying the condition $\Gamma(x,y)+\Gamma(y,x)=1$.

\subsection{Spectral properties}\label{ssec:IntroSpectral}
The main contribution of this paper lies in developing the spectral theory of digraphons. While some of our results follow from general principles in functional analysis, we have endeavored to present them in a manner accessible to researchers whose primary expertise is in graph theory.

Spectral graph theory traditionally studies the spectral properties of matrices associated with graphs, most notably the adjacency matrix.\footnote{Spectral graph theory also considers other matrices such as the Laplacian; however, these are not the focus of this work.} Specifically, given an adjacency matrix $A$ of a graph with $n$ vertices, one studies the properties of the linear operator $ A:\mathbb{C}^n \to \mathbb{C}^n $, defined for a row vector $\mathbf{v}$ by $ \mathbf{v} \mapsto \mathbf{v}A$. This framework extends naturally to the setting of graphons in a way we recall below.

Throughout the paper, we work within the complex Hilbert space $ L^2(\Omega) $, where $\Omega$ is a probability space equipped with a measure $\mu$ (and an implicit sigma-algebra). We write \index{$\Identity$ }$\Identity$ for the identity operator on a Banach space (the Banach space will be clear from the context). A digraphon $\Gamma$ defined on $\Omega$ induces an integral kernel operator $T_\Gamma: L^2(\Omega) \to L^2(\Omega) $, defined for $ f \in L^2(\Omega)$ as $g:=T_\Gamma f$ by
\begin{equation}\label{eq:defOper}
g(x) := \int_{\Omega} f(y)\Gamma(y, x) \, \mathrm{d}\mu(y), \quad \text{for all } x \in \Omega.
\end{equation}
\begin{defi}[eigenvalues, eigenfunctions, point spectrum, spectrum, spectral radius]\label{def:spectrum}
Suppose that $\Gamma$ is a digraphon on $\Omega$. A complex number $\tau$ and a function $f\in L^2(\Omega)$ are called an \emph{eigenvalue} and an \emph{eigenfunction} of $\Gamma$, respectively if $\Gamma(f)=\tau f$. The collection of all eigenvalues is called the \emph{point spectrum} of $\Gamma$, and denoted \index{$\PSpec(\Gamma)$}$\PSpec(\Gamma)$. The collection of all complex numbers $\tau$ for which $\Gamma-\tau\cdot\Identity$ is not invertible (as a bounded operator) is called the \emph{spectrum} of $\Gamma$, and denoted \index{$\Spec(\Gamma)$}$\Spec(\Gamma)$. The \emph{spectral radius} of $\Gamma$, denoted $\rho(\Gamma)$\index{$\rho(\cdot)$} is the supremum of moduli over all the eigenvalues of $\Gamma$ (and is defined as $\rho(\Gamma)=0$ if $\Spec(\Gamma)=\emptyset$).
\end{defi}
Definition~\ref{def:spectrum} is standard, see for example P7.3-5 in~\cite{MR0467220}. The following properties of $T_\Gamma$ are well known and will be used throughout without further reference.  We do not need precise definitions at this moment and refer to~\cite[Section~7.5]{MR3012035}.
\begin{fact}\label{fact:basicgraphon}
Suppose that $\Gamma$ is a digraphon. Then $T_\Gamma$ has the following properties:
\begin{enumerate}[label=(\roman*)]
    \item\label{en:BG1} $T_\Gamma$ is a compact Hilbert--Schmidt operator. 
    \item\label{en:BG2} We have $\Spec(\Gamma)\setminus\{0\}=\PSpec(\Gamma)\setminus\{0\}$. Consequently, the spectral radius may be defined by the point spectrum rather than the spectrum.
    \item\label{en:BG3} The spectrum is at most countable with the only possible accumulation point at~0,
    \item\label{en:BG4} Each eigenvalue has modulus at most $\|\Gamma\|_\infty$.    
    \item\label{en:BG5} If $\Gamma$ is a graphon then $T_\Gamma$ is self-adjoint. This in particular means that all the eigenvalues and eigenfunctions are real and that $\Gamma$ admits a \emph{spectral decomposition},
\begin{align}\label{eq:spectraldecomposition}
\Gamma(x, y) = \sum_{i} \lambda_i f_i(x) f_i(y)\;,
\end{align}
where $\{ \lambda_i \}_i$ are the eigenvalues of $ T_\Gamma $, and $ \{ f_i \}_i \subseteq L^2(\Omega) $ are the corresponding orthonormal eigenfunctions 
\end{enumerate}
\end{fact}
The operator-theoretic perspective on graphons was initiated in~\cite{borgs2012convergent} and further developed in~\cite{SZEGEDY20111156}. A comprehensive overview of the spectral theory of graphons is provided in Section~11.6 of~\cite{MR3012035}. In contrast, the spectral properties of digraphons have been investigated only recently, in~\cite{grebik2025convergencespectradigraphlimits}. In what follows, we summarize our main spectral results for digraphons.

Theorem~\ref{thm:periodicity} establishes a link between spectral properties and periodicity of a digraphon. The latter concept, defined in Definition~\ref{def:graphicallyperiodic}, is a natural counterpart to that known in finite digraphs and the theory of Markov chains; it concerns decomposing the ground space into `cyclic sets' so that all the arrows go only between cyclically consecutive cyclic sets. We prove that the periodicity of a strongly connected digraphon is equal to the number of peripheral eigenvalues (that is, eigenvalues whose modulus is the spectral radius, see Section~\ref{sssec:prelim_spectrum}).

In Proposition~\ref{prop:spectralradiusAndStrongComponents}, we prove that the spectral radius of a digraphon is equal to the maximum of spectral radii over all strong components.

Our third main result about spectral properties of concerns the asymptotic growth of the powers of a digraphon. Namely, in Theorem~\ref{thm:asymptotics} we prove that for a strongly connected aperiodic digraphon $\Gamma$ and $k$ large, for the $k$-th power of $\Gamma$ (see Definition~\ref{def:digraphonpower}) we have $\Gamma^k(x,y)\approx \rho(\Gamma)^k\cdot v_\textrm{left}(x)v_\textrm{right}(y)$, where $v_\textrm{left}$ and $v_\textrm{right}$ are the unique (see Theorem~\ref{thm:Schaefer74}\ref{en:skolka}) left and right eigenfunction of $\Gamma$ corresponding to eigenvalue $\rho(\Gamma)$.\footnote{The eigenfunction in Definition~\ref{def:spectrum} is the \emph{left} eigenfunction. The \emph{right} eigenfunction is the eigenfunction of the transposed digraphon $\Gamma^{\top}(x,y)=\Gamma(y,x)$. See Section~\ref{ssec:Duality} for more.}, and also give an extension to the periodic case. Theorem~\ref{thm:asymptotics} in fact applies to digraphons of nontrivial periodicity, too. This theorem is the most important ingredient needed in~\cite{HladkySavicky:Inhomo2SAT}, but we expect that it will have many additional applications. Indeed, as Definition~\ref{def:digraphonpower} tells us, $\Gamma^k$ expresses the homomorphism density of paths of length $k$, a quantity of obvious importance.

For the proof of Theorem~\ref{thm:asymptotics}, we introduce Proposition~\ref{prop:poweringPeriodic}. It asserts that if $\Gamma$ is strongly connected $D$-periodic, then the $D$ cyclic classes of $\Gamma$ are exactly strong components of $\Gamma^D$, and that $\Gamma^D$ is aperiodic. Further, it links the spectral properties of each component of $\Gamma^D$ to those of $\Gamma$.

\section{Preliminaries}
\subsection{Measure theory}
All subsets of measure spaces are assumed measurable, and all constructed subsets are evidently so. We use a measure space \index{$(\Omega,\mu)$}$\Omega$ equipped with a probability measure $\mu$ to this end. That is, the measure $\mu$ always implicitly underlies the space $\Omega$. For measurable subsets $A$ and $B$ in a measure space, we write \index{$=_0$, $\subseteq_0$}$A=_0 B$ and $A\subseteq_0 B$ for equality and containment modulo a nullset, respectively. 
Suppose that $f:\Omega\to\C$ is a function on a measure space $(\Omega,\mu)$. We write \index{$\supp f$}$\supp f$ for the \emph{support} of $f$, $\supp f=\{x\in \Omega: f(x)\neq 0\}$. When $f$ is real-valued, we write \index{$\essinf f$, $\esssup f$}$\essinf f$ and $\esssup f$ for the essential infimum and essential supremum of $f$, $\essinf f=\sup \{t\in \R:\mu(f^{-1}((-\infty,t))=0\}$, $\esssup f=\inf \{t\in \R:\mu(f^{-1}((t,+\infty))=0\}$. 

\subsection{Digraphons}\label{ssec:Digraphons}
We first introduce fundamental notions of the cut distance, digraphon representation of finite digraphs, subgraphons, degrees, and homomorphism densities within the digraphon framework. These definitions either exist in the literature or are direct analogues of their well-established counterparts in the symmetric (graphon) setting. For a comprehensive treatment of the graphon framework, see~\cite{MR3012035}.

First, we introduce the cut norm distance, the cut distance and a digraphon representation of a finite digraph.
Given two digraphons $U$ and $W$ on ground space $\Omega$, their \emph{cut norm distance} is defined as \index{$\cutnormdist(U,W)$}
\[\cutnormdist(U,W)=\sup_{S,T\subseteq \Omega}\left|\int_{S\times T} U-\int_{S\times T} W\right|\;.\] Note that in the case of graphons, often it is convenient to work with a different definition which involves only symmetric integrals $\int_{S\times S}$. It is well-known that for graphons the alternative definition differs from the original one by a factor of at most~4. This is however not true for digraphons, and only the non-symmetric definition is reasonable. To see this, consider two digraphons $A,B$ on ground set $[0,1]$, $A\equiv \frac12$ and $B(x,y)=\mathbbm{1}_{x\le y}$. Taking $S=[0,\frac12]$ and $T=[\frac12,1]$, we see that $\cutnormdist(U,W)\ge \frac18$ but $\sup_{S\subseteq \Omega}\left|\int_{S\times S} A-\int_{S\times S} B\right|=0$. The \emph{cut distance} is defined as \index{$\cutdist(U,W)$}
\begin{equation}\label{eq:defcutdistance}
\cutdist(U,W)=\inf_{\pi}\cutnormdist(U,W^\pi)\;,
\end{equation}
where $\pi$ ranges over measure preserving bijections on $\Omega$ and $W^\pi$ is a digraphon on $\Omega$ defined by $W^\pi(x,y)=W(\pi(x),\pi(y))$.

Next, we introduce digraphon representations of a finite digraph.
\begin{defi}\label{def:representationfinitedigraph}
Suppose that $G$ is a finite nonempty digraph. Partition the space $\Omega$ into $v(G)$ many sets of measure $\frac{1}{v(G)}$, $\Omega=\bigsqcup_{v\in V(G)}\Omega_v$. A \emph{digraphon representation} of $G$ is a digraphon on ground space $\Omega$ which is constant on each part $\Omega_v\times \Omega_w$, and whose value on that part is either~1 or~0, depending on whether $vw$ is a directed edge or not.
\end{defi}
A digraphon representation of a digraph is not unique as it depends on the choice of the partition $\Omega=\bigsqcup_{v\in V(G)}\Omega_v$. By the cut distance between a finite digraph $G$ and a digraphon $W$, we mean the cut distance of a digraphon representation of $G$ and $W$. Even though the digraphon representation is not unique as it depends on the choice of the partition $\Omega=\bigsqcup_{v\in V(G)}\Omega_v$, this does not change the cut distance due to the fact that~\eqref{eq:defcutdistance} involves measure preserving bijections.

We will use two ways of restricting a digraphon to subset of the ground set. 
\begin{defi}\label{def:restriction}
For a digraphon $U$ on $\Omega$ and for $A\subseteq \Omega$, let $U\llbracket A\rrbracket$\index{$U\llbracket A\rrbracket$} be a digraphon on $\Omega$ defined by 
\[
U\llbracket A\rrbracket(x,y)=\begin{cases}U(x,y) &\mbox{if $x,y\in A$, and}\\
0&\mbox{otherwise.}
\end{cases}
\]
Next, if $A$ has positive measure, we write $U_{\restriction A\times A}$\index{$U_{\restriction A\times A}$} for the restriction of $U$ to $A\times A$. We view $U_{\restriction A\times A}$ as a digraphon, that is, we assume that $A$ is a probability space equipped with the natural probability measure $\mu_A(\cdot)=\frac{\mu(\cdot\cap A)}{\mu(A)}$.
\end{defi}

Suppose $\Gamma$ is a digraphon on $\Omega$. For $x\in \Omega$, the \emph{indegree} and \emph{outdegree} of $x$ are defined as \index{$\degIn_\Gamma(x),\degIn_\Gamma(x)$}$\degIn_\Gamma(x):=\int_{y\in\Omega} \Gamma(y,x)\diff\mu(y)$ and $\degOut_\Gamma(x):=\int_{y\in\Omega} \Gamma(x,y)\diff\mu(y)$, respectively. There are variants of the indegree and outdegree \emph{relative to a given set $A\subseteq \Omega$}, $\degIn_\Gamma(x,A):=\int_{y\in A} \Gamma(y,x)\diff\mu(y)$ and $\degOut_\Gamma(x,A):=\int_{y\in A} \Gamma(x,y)\diff\mu(y)$. The \emph{minimum indegree}, \emph{minimum outdegree}, \emph{maximum indegree}, and \emph{maximum outdegree} of $\Gamma$ are defined by \index{$\degminIn(\Gamma)$, $\degmaxIn(\Gamma)$, $\degmaxIn(\Gamma)$, $\degmaxOut(\Gamma)$}$\degminIn(\Gamma):=\essinf_x \degIn_\Gamma(x)$ and $\degmaxIn(\Gamma):=\esssup_x \degIn_\Gamma(x)$, and similarly for $\degminOut(\Gamma)$ and $\degmaxOut(\Gamma)$.

Suppose that $D$ is an oriented graph on vertex set $[n]$. The \emph{homomorphism density of $D$ in $\Gamma$} is defined as
\begin{equation}\label{eq:defHomDens} 
\index{$t(D,\Gamma)$}
t(D,\Gamma)=\int_{x_1}\int_{x_2}\ldots\int_{x_n}\prod_{(i,j) \in E(D)}\Gamma(x_i,x_j)\diff\mu^n\;.
\end{equation}

The Counting lemma for graphons asserts that graphons that are close in the cut distance have similar homomorphism densities of smallish graphs. Here, we include a digraphon counterpart. We omit the proof, since it follows by the same telescoping argument (see, e.g., Lemma~10.23 in~\cite{MR3012035}). However, as we remark below, the version of the counting lemma holds only for oriented graphs, that is, digraphs without counterparallel edges.
\begin{lem}\label{lem:CountingLemma}
    Suppose that $U$ and $W$ are two digraphons, and $D$ is an oriented graph. Then $|t(D,U)-t(D,W)|\le e(D)\cdot \cutdist(U,W)$.
\end{lem}
\begin{rem}
Consider $n$ large. Let $G_1$ be a random orientation of the complete graph $K_n$. Let $G_2$ be a digraph created from the Erdős--Rényi random graph with edge probability~$\frac12$ by replacing every edge with directed edges in both directions. In is easy to show that digraphon representations $W_1$ and $W_2$ of $G_1$ and $G_2$ have the cut distance $o(1)$ asymptotically almost surely. Yet, obviously for the 2-cycle we have $t(C_2,W_1)=0$ and $t(C_2,W_2)=\frac12+o(1)$ with high probability. This example shows that the restriction to oriented graphs in Lemma~\ref{lem:CountingLemma} is necessary.
\end{rem}

Suppose that $k\in \N$. Write \index{$P_k,P_k^{\bullet\bullet}$}$P_k$ for the directed path $1,2,3,\cdots ,k,{k+1}$, and $P_k^{\bullet\bullet}$ for the directed path rooted at its two terminal vertices. Write \index{$C_k,C_k^\bullet$}$C_k$ for the directed cycle $1,2,3,\cdots ,k$ and $C_k^\bullet$ for $C_k$ rooted at vertex $1$. We define \index{$t^{\bullet\bullet}_{x,y}(P_k^{\bullet\bullet},\Gamma)$}$t^{\bullet\bullet}_{x_1,x_{k+1}}(P_k^{\bullet\bullet},\Gamma)$ as a function of $x_1,x_{k+1}\in\Omega$ by disintegrating $t(P_k,\Gamma)$ with respect to $x_1$ and $x_{k+1}$. Likewise, we define \index{$t^{\bullet}_{x}(C_k^{\bullet},\Gamma)$}$t^{\bullet}_{x_1}(C_k^{\bullet},\Gamma)$ as a function of $x_1\in\Omega$ by disintegrating $t(C_k,\Gamma)$ with respect to $x_1$. That is,
\begin{align}    
\label{eq:densPathTwoTerminals}
t^{\bullet\bullet}_{x_1,x_{k+1}}(P_k^{\bullet\bullet},\Gamma)&=\int_{x_2}\int_{x_3}\ldots\int_{x_k}\prod_{i=1}^{k}\Gamma(x_i,x_{i+1})\;\diff\mu^{k-1}\;,\\
\nonumber
t^{\bullet}_{x_1}(C_k^{\bullet},\Gamma)&=\int_{x_2}\int_{x_3}\ldots\int_{x_k}\prod_{i=1}^{k-1}\Gamma(x_i,x_{i+1})\cdot \Gamma(x_k,x_1)\;\diff\mu^{k-1}\;.
\end{align}

Last, we introduce the power of a digraphon. This definition is similar to the definition of matrix powers.
\begin{defi}[Power of a digraphon]\label{def:digraphonpower}
Suppose that $\Gamma$ is a digraphon on $\Omega$, and let $k\in\N$. Define \index{$\Gamma^k$}$\Gamma^k$ as digraphon on $\Omega$ by $\Gamma^k(x,y):=t^{\bullet\bullet}_{x,y}(P_{k}^{\bullet\bullet},\Gamma)$.
\end{defi}
It is straightforward to verify that the definition is consistent with operator powers, that is, for every $f\in L^2(\Omega)$ and every $k\in \N$ we have
\begin{equation}\label{eq:powerfordummies}
\Gamma^k f=\underbrace{\Gamma\Gamma\cdots \Gamma}_{\mbox{$k$ times}} f\;.
\end{equation}
We highlight an important notational difference: $\Gamma^k(x,y)$ refers to Definition~\ref{def:digraphonpower} whereas $\Gamma(x,y)^k$ means taking the real number $\Gamma(x,y)$ and rising it to the $k$-th power. When written without parameters, $\Gamma^k$ always refers to Definition~\ref{def:digraphonpower}.

\subsubsection{Calculus of reachability}
A central question in the study of finite digraphs is whether a directed path exists between two vertices, and if so, how many. In this section, we introduce tools for investigating the analogous problem in the setting of digraphons.

Suppose that $\Gamma$ is a digraphon on $\Omega$. Let $x,y\in \Omega$. We say that \emph{$x$ is reachable from $y$} if for some $k\in\N$ we have $\Gamma^k(y,x)>0$. The \emph{reachability sequence from $x$ to $y$} is the set of all such numbers $k$.
We write 
\index{$\ReachOut_\Gamma(y),\ReachIn_\Gamma(y)$,$\ReachOut_\Gamma(Y),\ReachIn_\Gamma(Y)$}
$\ReachOut_\Gamma(y)\subseteq \Omega$ for all $x$'s that are reachable from $y$. We write $\ReachIn_\Gamma(x)\subseteq \Omega$ for all $y$'s so that $x\in\ReachOut_\Gamma(y)$. For a set $Y$ we write $\ReachOut_\Gamma(Y)\subseteq \Omega$ for the set of all $x$'s such that $\mu(\{y\in Y:x\in \ReachOut_\Gamma(y)\})>0$. The set $\ReachIn_\Gamma(Y)$ is defined similarly.

The next lemma is a natural transitivity principle for reachability.

\begin{lem}\label{lem:calculusofreachability}
Suppose that $\Gamma$ is a digraphon on $\Omega$, $x\in \Omega$, and $Y\subseteq \Omega$, and $c\in[0,1)$. Define a function $h:\Omega\to[0,1]$ by $h(z):=\mu(\ReachOut_\Gamma(z)\cap Y)$.
Suppose that for a positive measure of $y$'s that are reachable from $x$ we have that $h(y)>c$. Then we have $h(x)>c$.
\end{lem}
\begin{proof}
Let $T\subseteq \ReachOut_\Gamma(x)$ be the set of $y$'s with $h(y)>c$. By the assumption, $\mu(T)>0$. Let 
$$
S:=\{s\in Y\::\: \mu(T\cap \ReachIn(s))>0\}
\;.
$$
Let us use Fubini's Theorem to express the measure $\mu^2(\{(y,s)\in T\times S: s\in \ReachOut_\Gamma(y)\})$. On the one hand, for almost every $y\in T$, we have $\mu(\{s\in S:s\in \ReachOut_\Gamma(y)\})>c$. On the other hand, for each $s\in S$, we trivially have $\mu(\{y\in T:s\in \ReachOut_\Gamma(y)\})\le \mu(T)$. Hence,
\[
\mu(T)c
<
\mu^2(\{(y,s)\in T\times S: s\in \ReachOut_\Gamma(y)\})
\le
\mu(T)\mu(S)\;,
\]
which gives $\mu(S)>c$. It is easy to see that $S\subseteq_0 \ReachOut_\Gamma(x)$. The claim follows.
\end{proof}

\subsubsection{Regularity Lemma}
Szemerédi's Regularity Lemma has played pivotal roles in graph theory for several decades. Slight variants tailored for different purposes exist. Here, we include a minimum outdegree variant of the Weak Regularity Lemma for digraphons. Its proof follows the same line as the original proof.

\begin{defi}\label{def:clusterdigraph}
Suppose that $\Gamma$ is a digraphon on $\Omega$. Suppose that $\Omega=Z_1\sqcup Z_2\sqcup \ldots\sqcup Z_t$ is a partition into sets of measures $\frac1t$. Suppose that $d,\eps>0$ are given. We say that a digraph $D$ on vertex $[t]$ is a \emph{$(d,\eps)$-cluster digraph} for partition $\Omega=Z_1\sqcup Z_2\sqcup \ldots\sqcup Z_t$ and digraphon $\Gamma$, if the following conditions hold:
\begin{enumerate}[label=(\roman*)]
    \item For each edge $ij$ of $D$, the average value $p_{i,j}:=t^2\int_{Z_i\times Z_j} \Gamma$ of $\Gamma$ on $Z_i\times Z_j$ is at least $d$.
    \item\label{en:clCut} For each edge $ij$ of $D$, consider two digraphons $\Gamma_{i,j}$ and $\Gamma^{avg}_{i,j}$ on $\Omega$, defined as follows: We have $\Gamma_{i,j}(x,y)=\Gamma^{avg}_{i,j}(x,y):=0$ for $(x,y)\not\in Z_i\times Z_j$. For $(x,y)\in Z_i\times Z_j$, we define $\Gamma_{i,j}(x,y):=\Gamma(x,y)$ and $\Gamma^{avg}_{i,j}(x,y):=p_{i,j}$. Then we have that $\cutnormdist(\Gamma_{i,j},\Gamma^{avg}_{i,j})<\eps\cdot t^{-2}$.
    \item The minimum outdegree of $D$ is at least $(\degminOut(\Gamma)-d-\eps)t$.
\end{enumerate}
\end{defi}
The regularity lemma asserts that digraphons have cluster digraphs of bounded order.
\begin{thm}\label{thm:Regularity}
    Suppose that $\eps>0$ is given. Then there exists a number $t_0\in \N$ so that for every $d\in [0,1]$ and for every digraphon on $\Omega$ there exists a partition $\Omega=Z_1\sqcup Z_2\sqcup \ldots\sqcup Z_t$ into sets of measures $\frac1t$ for some $t\le t_0$ and a $(d,\eps)$-cluster digraph for this digraphon and partition.
\end{thm}

\subsection{Spectral theory}\label{ssec:SpectralTheory}
As explained in Section~\ref{ssec:IntroSpectral}, each digraphon can be viewed as an integral kernel operator via~\eqref{eq:defOper}. In this section, we introduce spectral tools tailored to the study of digraphons. To the best of our knowledge, spectral methods for digraphons had previously been only used in~\cite{GRZESIK2023117,InhomogeneousRandomDigraphs}, where their use was rather elementary, and recently, in a focused way in~\cite{grebik2025convergencespectradigraphlimits}, which proves continuity of the spectra of digraphons with respect to the cut distance.
In contrast, our approach relies on more sophisticated techniques drawn from the theory of Banach lattices. While the results presented in this section are mostly borrowed from the literature, assembling and adapting them to the setting of digraphons has been a nontrivial task --- especially for us, as non-specialists in functional analysis. In fact, translating abstract functional analytic language into a form applicable to digraphons proved to be one of the most challenging aspects of this work.
It is somewhat paradoxical that this translation process is circular in nature. That is, the development of Hilbert space theory in the early 20th century, as well as the emergence of Banach lattice theory in the 1940s and 1950s, was originally motivated by the study of integral kernel operators. Yet, most modern treatments of these topics adopt an abstract perspective, requiring us to reinterpret them in the concrete setting of digraphons (which are themselves integral kernel operators).
Throughout this section, we aim to make the relevant functional analytic tools as accessible as possible to readers who, like us, may not be specialists in the area.

As mentioned earlier, all our Hilbert and Banach spaces are complex. This reflects the well-known fact that in the non-self-adjoint setting (which necessarily comes with non-symmetric digraphons), real Hilbert spaces lack many favorable features. We write \index{$\ImaginaryUnit$}$\ImaginaryUnit$ for the imaginary unit.

\subsubsection{Spectrum}\label{sssec:prelim_spectrum}
Definition~\ref{def:spectrum} extends verbatim for every bounded linear operator $T: L^2(\Omega) \to L^2(\Omega)$. While not all operators in this paper will be digraphons, all of them will be compact. In particular, properties stated as Fact~\ref{fact:basicgraphon}\ref{en:BG2} and~\ref{en:BG3} hold for such an operator, too.

Let us add three more spectral definitions. An eigenvalue $\lambda$ and its corresponding eigenfunction $f$ are \emph{peripheral} if $|\lambda|=\rho(T)$.
If $\lambda$ is an eigenvalue, then its \emph{algebraic multiplicity} is
\index{$\mathfrak{m}_T(\lambda)$}
\begin{equation}
\mathfrak{m}_T(\lambda):=\dim\left( \bigcup_{n=1}^\infty \ker\left((T - \lambda\cdot \Identity)^n\right) \right)
\end{equation}
We say that an eigenvalue $\lambda$ is \emph{simple} if its algebraic multiplicity equals~1.

We shall use one of the main results of a recent paper of Grebík, Král', Liu, Pikhurko, and Slipantschuk~\cite{grebik2025convergencespectradigraphlimits}, which expresses the homomorphism density of a cycle in terms of its spectrum. Let us note that for graphons, the same formula has been known (see Equation~(7.23) in~\cite{MR3012035}).
\begin{prop}\label{prop:PU}
For any digraphon $\Gamma$ and $k\ge 3$ we have 
\begin{equation*}
t(C_k,\Gamma)=\sum_{\lambda\in\PSpec(\Gamma)}\mathfrak{m}_{\Gamma}(\lambda)\lambda^k
\;.
\end{equation*}
\end{prop}

\subsubsection{Kernels}
A \emph{kernel} is an arbitrary (complex) function $U\in L^\infty(\Omega)$. We will use kernels in the proof of Theorem~\ref{thm:asymptotics}; those kernels will be of the form $U=W_1-W_2$, where $W_1$ and $W_2$ are digraphons. The definitions of homomorphism densities~\eqref{eq:defHomDens}, ~\eqref{eq:densPathTwoTerminals}, as well as of powers (Definition~\ref{def:digraphonpower} and property~\eqref{eq:powerfordummies}) extend to kernels.

\subsubsection{Duality}\label{ssec:Duality}
Suppose that $X$ is a Banach space. Recall that $X^*$ denotes the \emph{dual of $X$}, which is a Banach space comprising of all continuous linear functionals $X\to \C$. If $T:X\to X$ is a bounded linear operator on $X$, then the \emph{dual operator} $T^*:X^*\to X^*$ is defined by
\[
(T^*\psi)(x)=\psi(Tx) \qquad \text{for all }\psi\in X^*,\; x\in X.
\]

The relevant setting for us will be when $X=L^2(\Omega)$. In that case, the Riesz Representation Theorem asserts that for each continuous functional $\psi\in X^*$ there exists its \emph{Riesz representation} $g_\psi\in L^2(\Omega)$ so that
\[
\psi(f)=\int f(y)\,\overline{g_\psi(y)}\,\diff\mu(y) \qquad \text{for every } f\in L^2(\Omega).
\]
This allows us to identify $(L^2(\Omega))^*=L^2(\Omega)$.

Under this identification, the dual operator $T^*$ coincides with the Hilbert space adjoint in the usual sense: for all $f,g\in L^2(\Omega)$ one has
\begin{equation}\label{eq:Adjoint}    
\langle Tf , g\rangle = \langle f, T^* g\rangle .
\end{equation}
This identity follows directly from the definition of $T^*$ on functionals together with the Riesz representation.

If $\Gamma$ is a digraphon viewed as an integral kernel operator $T_{\Gamma}:L^2(\Omega)\to L^2(\Omega)$, then the adjoint operator $T_\Gamma^*:L^2(\Omega)\to L^2(\Omega)$ is the integral kernel operator for the transposed digraphon $\Gamma^{\top}$, defined by
\[
\Gamma^{\top}(x,y)=\Gamma(y,x).
\]
Eigenfunctions and eigenvalues of $\Gamma^{\top}$ are sometimes called \emph{right eigenfunctions} and \emph{right eigenvalues} of $\Gamma$. To emphasize the difference, the original eigenfunctions and eigenvalues of $\Gamma$ may be referred to as \emph{left}.

\subsubsection{Gelfand's formula for the spectral radius}
We use the famous Gelfand's formula  (see P7.5-5 in~\cite{MR0467220}). Recall that the \emph{operator norm} $\|R\|_{\mathrm{op}}$ of an operator $R$ on a Banach space with norm $\|\cdot\|$ is defined $\|R\|_{\mathrm{op}}:=\sup_{x:\|x\|\le 1}\|Rx\|$.
\begin{prop}[Gelfand's formula for the operator norm]\label{pro:Gelfand}
Suppose that $T$ is a bounded operator on a complex Banach space. Then we have
\begin{equation}\label{eq:spectralradius}
    \rho(T)=\lim_{k\to\infty}\sqrt[k]{\|T^k\|_{\mathrm{op}}}\;.
\end{equation}    
\end{prop}
Gelfand's formula~\eqref{eq:spectralradius} holds not only for the operator norm, but also for many other norms. The version we need here is for the $L^2$-norm.
\begin{prop}[Gelfand's formula for the $L^2$-norm]\label{pro:GelfandHS}
Suppose that $\Gamma\in L^2(\Omega^2)$ is given. Then,
\begin{equation}\label{eq:gelfand-1}
    \rho(\Gamma)=\lim_{k\to\infty}\sqrt[k]{\|\Gamma^k\|_{2}}\;.
\end{equation}    
\end{prop}
This statement may be standard, but we were not able to find it in literature. We include a proof communicated to us by Vladimír Müller.
\begin{proof}
In the proof, we are going to work with the Hilbert--Schmidt norm $\|\cdot\|_{\mathrm{HS}}$, which is defined for some bounded operators on $L^2(\Omega)$. As we are concerned with digraphons, we shall need the Hilbert--Schmidt norm only in the case of integral kernel operators, and thus we introduce it only in this particular setting. All the properties of the Hilbert--Schmidt norm we shall need can be found on page~267 of~\cite{MR1070713}. In particular, the Hilbert--Schmidt norm is always (when it is defined) an upper-bound on the operator norm, $\|\cdot \|_{\mathrm{op}}\le \|\cdot \|_{\mathrm{HS}}$. Also, for integral kernel operators, the Hilbert--Schmidt norm is equal to the $L^2$-norm of the corresponding kernel. Combining this with~\eqref{eq:spectralradius}, we have $\rho(\Gamma)\le \liminf_{k\to\infty}\sqrt[k]{\|\Gamma^k\|_{\mathrm{HS}}}=\liminf_{k\to\infty}\sqrt[k]{\|\Gamma^k\|_{2}}$. So, it remains to establish the other inequality. The last fact from page~267 of~\cite{MR1070713} we use is that $\|A \cdot B\|_{\mathrm{HS}}\le \|A \|_{\mathrm{op}}\cdot \|B \|_{\mathrm{HS}}$. Applying this with $A=T^{k-1}$ and $B=T$, we have
\[
\limsup_{k\to\infty}\sqrt[k]{\|\Gamma^k\|_{\mathrm{HS}}}
\le
\limsup_{k\to\infty}\sqrt[k]{\|\Gamma^{k-1}\|_{\mathrm{op}}}
\cdot 
\lim_{k\to\infty}\sqrt[k]{\|\Gamma\|_{\mathrm{HS}}}\;.
\]
The first term tends to $\rho(\Gamma)$ by~\eqref{eq:spectralradius}. The second term tends to~1. Hence, the proof follows.
\end{proof}

We shall use the well-known fact that the spectral radii of a bounded operator and its dual are the same, see, e.g., Proposition~6.1 in~\cite{MR1070713}.
\begin{prop}\label{prop:SpectralRadiusDual}
If $T$ is a bounded operator on a Hilbert space and $T^*$ is its dual, then $\rho(T)=\rho(T^*)$.
\end{prop}

The next lemma concerns $L^\infty$-boundedness of eigenfunctions of digraphons.
\begin{lem}\label{lem:eigenvectorsbounded}
Suppose that $\Gamma$ is a digraphon on $\Omega$, and $f$ an eigenfunction for eigenvalue $\gamma\neq0$. Then $\|f\|_\infty\le \frac{\|\Gamma\|_\infty \|f\|_1}{|\gamma|}$.
\end{lem}
\begin{proof}
We have $\gamma f=\Gamma f$. Fix an arbitrary $x\in\Omega$. Thus, $|\gamma f(x)|=|\int_y f(y)\Gamma(y,x)|\le \|\Gamma\|_\infty\cdot\|f\|_1$.
\end{proof}
Our last result of this section concerns powers (in the sense of Definition~\ref{def:digraphonpower}) of a kernel. Namely, it asserts that the $L^\infty$-norm of a high power of a kernel does not grow much faster than the corresponding power of its spectral radius.
\begin{prop} \label{prop:asymptotics}
Suppose that $U\in L^\infty(\Omega^2)$ is a kernel, and $\beta>\rho(U)$. Then, as $k\to\infty$, we have $\| U^k \|_\infty = O(\beta^k)$.
\end{prop}
\begin{proof}
For each $(x,y)\in \Omega$ and $k\ge 2$, we have,
$$
U^k(x,y) = \int_{\Omega^2} U(x,z) U^{k-2}(z,t) U(t,y) \diff z \diff t\;,
$$
and therefore
$$
|U^k(x,y)| \le \int_{\Omega^2} |U(x,z)| \; |U^{k-2}(z,t)| \; |U(t,y)| \diff z \diff t \;.
$$
Since $U$ is bounded, we have uniformly over $(x,y)$ and $k$,
\begin{align*}
|U^k(x,y)| &= O\left( \int_{\Omega^2} |U^{k-2}(z,t)| \diff z \diff t \right)
= O\left( \| U^{k-2} \|_1 \right)\\
\JUSTIFY{$\|\cdot\|_1\le \|\cdot\|_2$  on a probability space}&= O\left( \| U^{k-2} \|_2 \right)\\
\JUSTIFY{by~\eqref{eq:gelfand-1}, $k$ large}&=O(\beta^k)\;,
\end{align*}
as was needed.
\end{proof}

\subsubsection{Banach lattices}\label{ssec:BanachLattices}
Banach lattices are partially ordered Banach spaces. We do not need to recall the general axioms here, as we work exclusively with the Banach space $L^2(\Omega)$ in which the order is given by $f\le g$ if $f(x)\le g(x)$ for every $x\in \Omega$.
The second lemma connects the notion of strongly connected graphons with the concept of irreducible operators in the theory of Banach lattices. We give necessary definitions from \cite{MR423039}, but only in the specific setting of the Banach lattice $L^2(\Omega)$. For $f\in L^2(\Omega)$, we write $|f|$ for the pointwise absolute value. A set $I\subseteq L^2(\Omega)$ is an \emph{ideal} if it is a closed subspace of $L^2(\Omega)$ such that for every $f\in I$ and for every $g\in L^2(\Omega)$ with $|g|\le |f|$ we have $g\in I$. Fortunately, ideals in $L^2(\Omega)$ have a simple structure. Namely (see, e.g., \cite[Theorem~7.10]{MR3410920}), every ideal is of the form $I_Z=\{f\in L^2(\Omega):f_{\restriction Z}\equiv 0\}$ for some measurable $Z\subseteq \Omega$. An operator $T$ on $L^2(\Omega)$ is \emph{irreducible} if the only ideals $I$ for which $TI\subseteq I$ are $I=\{0\}$ and $I=L^2(\Omega)$. Below, we connect this concept to digraphons by proving that strongly connected digraphons are irreducible. The converse also holds and is straightforward to show, but we will not need it here.
\begin{lem}\label{lem:irreducible}
Suppose that $\Gamma$ is a strongly connected digraphon on $\Omega$. Then $\Gamma$ is irreducible.
\end{lem}
\begin{proof}
Let $I$ be an ideal. Assume at $I$ is of the form $I_Z$ as above. There is nothing to prove if $\mu(Z)\in\{0,1\}$, so assume $\mu(Z)\in (0,1)$. By the definition of strong connectedness, we have $\int_{(\Omega\setminus Z)\times Z}\Gamma>0$. That means that for the characteristic function $\mathbbm{1}_{\Omega\setminus Z}\in I_Z$ we have that the support of $\Gamma \mathbbm{1}_{\Omega\setminus Z}$ nontrivially intersects $Z$. Thus, $\Gamma\mathbbm{1}_{\Omega\setminus Z}\notin I_Z$. We have concluded that the only ideals with $\Gamma I\subseteq I$ are the two trivial ones.
\end{proof}

\subsubsection{Reducing pairs}\label{ssec:ReducingPairs}
In this section, we deal with decompositions of operators. Let us recall some basic notation. If $X$ is a Banach space and $X_1$ and $X_2$ are closed linear subspaces of $X$, we say that $X$ is a  \emph{direct sum} of $X_1$ and $X_2$, written as $X=X_1\oplus X_2$, if every $x\in X$ can be written in a unique way as $x=x_1+x_2$, with $x_i\in X_i$. A bounded linear operator $P$ on $X$ is a \emph{projection} if $P^2=P$. Suppose that $Y$ is a subspace of $X$. We say that a projection $P$ is \emph{on $Y$} if $\image(P)=Y$. The next fact is well-known.
\begin{fact}\label{fact:projectionsId}
Suppose that $P_1$ and $P_2$ are the canonical projections associated with the decomposition
$X = X_1 \oplus X_2$.
\end{fact}
\begin{proof}
Since $X = X_1 \oplus X_2$, every $x \in X$ has a unique decomposition $x = x_1 + x_2$, $x_1 \in X_1$ and $x_2 \in X_2$.
We have $P_1(x) = x_1$ and $P_2(x) = x_2$. Therefore,
\[
(P_1 + P_2)(x) = P_1(x) + P_2(x) = x_1 + x_2 = x\;.
\]
Since $x$ was arbitrary, we get the statement.
\end{proof}

Suppose that $T$ is a bounded linear operator on a Banach space $X$. Suppose that we can write $X$ as a direct sum of closed subspaces, $X=X_1\oplus X_2$. We say that $(X_1,X_2)$ is a \emph{reducing pair for $T$} if $TX_1\subseteq X_1$ and $TX_2\subseteq X_2$. This definition is standard, see for example page~81 in~\cite{MR1921782}. The proposition follows from Theorem~2.22 in~\cite{MR1921782}. 
\begin{prop}\label{prop:reducingpair}
Suppose that $(X_1,X_2)$ is a reducing pair for a bounded linear operator $T$ on a Banach space $X$. Then there exist projections $P_1$ and $P_2$ on $X$ such that for $i=1,2$ we have $\image(P_i)=X_i$ and $T\circ P_i=P_i \circ T$.
\end{prop}
In linear algebra, reducing pairs arise naturally: whenever a linear operator can 
be written in block diagonal form, the corresponding decomposition of the space 
into the two invariant blocks provides a reducing pair.  This familiar situation 
motivates the use of reducing pairs in the functional-analytic setting, where 
they play an analogous role in decompositions associated with spectral sets.

Suppose that $T$ is a bounded linear operator on $X$. A set $\sigma \subseteq \Spec(T)$ is called a \emph{spectral set} if it is closed and open in the relative topology on $\Spec(T)$. This definition is standard, see for example Definition~6.32 in~\cite{MR1921782}. We then have the following fundamental result, see Theorem~6.34 in~\cite{MR1921782}.
\begin{thm}\label{thm:reducing}
Suppose that $T$ is a bounded linear operator on $X$. Suppose that $\sigma\notin\{ \emptyset,\Spec(T)\}$ is a spectral set of $T$. Then $T$ admits a reducing pair $(X_1,X_2)$ into two nonempty subspaces such that $\Spec(T_{\restriction X_1})=\sigma$ and $\Spec(T_{\restriction X_2})=\Spec(T)\setminus \sigma$.
\end{thm}

We give a short description of the structure of a reducing pair in the special case when the spectral set is a singleton of algebraic multiplicity~1 (which is what we need).
\begin{lem}\label{lem:ReducingPairAM}
Suppose that $T$ is a compact bounded linear operator on $X$. Suppose that $\sigma=\{\lambda\}$ is a spectral set of $T$, where $\lambda\neq 0$ has algebraic multiplicity~1. Suppose that $(X_1,X_2)$ is the corresponding reducing pair. Then $X_1=\ker(\lambda\Identity -T)$, $X_2=\image(\lambda\Identity -T)$, and $\dim X_1=1$.
\end{lem}
\begin{proof}
The proof relies on functional calculus and we only give a sketch, referring to~\cite{MR1921782} for details, including definitions of the resolvent function, and of the pole. Firstly, Theorem~6.39 in~\cite{MR1921782} tells us that $\lambda$ is a pole of order~1 of the resolvent function of~$T$. Corollary~6.40 in~\cite{MR1921782} then tells us that $X_1=\ker(\lambda\Identity -T)$ and $X_2=\image(\lambda\Identity -T)$. As the algebraic multiplicity of $\lambda$ is~1, so is its geometric multiplicity. In particular, $\dim(\ker(\lambda\Identity -T))=1$.
\end{proof}

\subsubsection{The \Krein--Rutman theorem and beyond}
The Perron--Frobenius Theorem asserts that a real square matrix with nonnegative entries has a unique non-negative eigenvector, and that this eigenvector corresponds to the eigenvalue of the largest eigenvalue in absolute value. The \Krein--Rutman Theorem is often considered a counterpart of the Perron--Frobenius Theorem for nonnegative operators. For main results of this paper regarding periodicity of digraphons (Theorem~\ref{thm:periodicity}) and asymptotic growth (Theorem~\ref{thm:asymptotics}) we need a somewhat strengthened version of the \Krein--Rutman Theorem which we give in Theorem~\ref{thm:Schaefer74} below. While we could not find Theorem~\ref{thm:Schaefer74} in this form in literature, it is just a slight variation of known results, as we show later. Given a digraphon $\Gamma$ on $\Omega$, write \index{$\ground(\Gamma)$}$\ground(\Gamma)=\{x\in\Omega:\int_\Omega \Gamma(x,y)\diff \mu(y)>0\}\cup \{y\in\Omega:\int_\Omega \Gamma(x,y)\diff \mu(x)>0\}$. 
\begin{thm}\label{thm:Schaefer74}
Suppose that $\Gamma$ is a digraphon on $\Omega$ with $\rho(\Gamma)>0$ and such that $\ground(\Gamma)$ is strongly connected.
Then for some $d \in\N$, we have the following.
\begin{enumerate}[label=(\roman*)]
\item\label{en:EigenvaluesAndSimple} The set of eigenvalues of modulus $\rho(\Gamma)$ is
$\left\{\exp(-2\pi \ImaginaryUnit k/d) \rho(\Gamma)\::\:k=0,\ldots,d-1\right\}$. Each of these eigenvalues is simple.
\item\label{en:LeadingEigenNonnegative} There are right and left non-negative eigenfunctions $v_\mathrm{right}$ and $v_\mathrm{left}$ for the eigenvalue $\rho(\Gamma)$. 
\item\label{en:LeadingEigenvectorPositive} We have $\supp(v_\mathrm{right})=_0\supp(v_\mathrm{left})=_0\ground(\Gamma)$.
\item\label{en:skolka} The eigenvalue $\rho(\Gamma)$ is the only one with a real and non-negative eigenfunction.
\item\label{en:PrincipalEigevectorInequality} If $f\in L^2(\Omega)$ is nonnegative such that $\Gamma f\ge \rho(\Gamma)f$ then $f$ is equal to $v_\mathrm{right}$ (up to constant multiple).
\end{enumerate}
\end{thm}

Two results that will be used to derive Theorem~\ref{thm:Schaefer74}. The first one is the original \Krein--Rutman Theorem.
\begin{thm}\label{thm:KreinRutman}
    Suppose that $T$ is a compact operator on $L^2(\Omega)$ which is nonnegative. Suppose that $T$ has a non-zero eigenvalue. Then there exists a nonnegative $f\in L^2(\Omega)$ such that $Tf=\rho(T)f$.
\end{thm}
The second one is Theorem~V.5.2 in~\cite{MR423039}. We rephrase this result in Theorem~\ref{thm:ActualSchaefer} below in our language, expanding necessary definitions. We also narrow the statement down to the Banach space $L^2(\Omega)$ even though the original statement is for a general Banach lattice, and to integral kernel operators. Since our integral kernel operators come from digraphons (which are nonnegative), they are automatically `positive' in the sense of Definition~II.2.4 in~\cite{MR423039}. This is a basic assumption for the whole theory of Banach lattices in~\cite{MR423039}. 
\begin{thm}\label{thm:ActualSchaefer}
Suppose that $\Gamma$ is a digraphon on $\Omega$ of spectral radius $\rho(\Gamma)>0$. Suppose that $\Gamma$ is an irreducible operator. Suppose that the dual operator $\Gamma^*$ of $\Gamma$ possesses an eigenfunction $\psi\in (L^2(\Omega))^*$ for eigenvalue~$\rho(\Gamma)$ such that its Riesz representation $g_\psi$ is a nonnegative function.\footnote{This part is a translation of `possesses an invariant form' in Theorem~V.5.2 in~\cite{MR423039} to our setting}. The following assertions are true:
\begin{enumerate}[label=(\roman*)]
    \item The eigenvalues of $\Gamma$ of modulus~$\rho(\Gamma)$ are a subgroup of the circle group.
    \item\label{en:rotationsymmetry} Each eigenvalue $\gamma$ of $\Gamma$ of modulus~$\rho(\Gamma)$ is simple. Further, $\PSpec(\Gamma)=\gamma\cdot\PSpec(\Gamma)$ including algebraic multiplicities.
    \item The only eigenvalue of $\Gamma$ with a nonnegative eigenfunction is~$\rho(\Gamma)$.
\end{enumerate}
\end{thm}
\begin{proof}[Proof of Theorem~\ref{thm:Schaefer74}]
First we apply Theorem~\ref{thm:KreinRutman} on $\Gamma$ to find its eigenfunction $v_\mathrm{right}$ for the eigenvalue $\rho(\Gamma)$. The theorem also asserts that $v_R$ is non-negative. We then switch to the dual operator $\Gamma^*$. Recall that by Proposition~\ref{prop:SpectralRadiusDual}, its spectral radius is the same. We again find its eigenfunction $v_\mathrm{left}$ for eigenvalue $\rho(\Gamma)$. By the discussion above this is the left eigenvalue of $\Gamma$.

We now prove that $\supp(v_\mathrm{right})=\ground(\Gamma)$. First, let us prove that $\supp(v_\mathrm{right})\subseteq \ground(\Gamma)$. Indeed, if $x\notin \ground(\Gamma)$, then in the identity $\rho(\Gamma)f(x)=\int_y \Gamma(x,y)f(y)$, the term $\Gamma(x,y)$ is~0 almost every. Hence, $f(x)=0$. Now, suppose that $A:=\supp(v_\mathrm{right})$ has measure strictly less than (on the other hand, obviously, the measure of $A$ is positive) $\mu(\ground(\Gamma))$. Take $B:=\ground(\Gamma)\setminus A$. Since by strong connectedness, $\int_{A\times B}\Gamma(x,y)>0$, we also have $\int_{A\times B}\Gamma(x,y)f(x)>0$. This means that $\int_B f(y)>0$, a contradiction.

Part~\ref{en:LeadingEigenvectorPositive} follows from Theorem~\ref{thm:ActualSchaefer}. Last, we turn to Part~\ref{en:PrincipalEigevectorInequality}, which is also standard. Set $g:=\Gamma f- \rho(\Gamma)f$. We have $g\ge 0$. Consider the inner product
\begin{align*}
\langle g,v_\mathrm{left}\rangle
&=
\langle T_\Gamma f- \rho(\Gamma)f,v_\mathrm{left}\rangle
=
\langle T_\Gamma f,v_\mathrm{left}\rangle
-
\langle \rho(\Gamma)f,v_\mathrm{left}\rangle
=
\langle f,T_\Gamma^* v_\mathrm{left}\rangle
-
\langle \rho(\Gamma)f,v_\mathrm{left}\rangle
\\
&=
\rho(\Gamma)\langle f,v_\mathrm{left}\rangle
-
\rho(\Gamma)\langle f,v_\mathrm{left}\rangle
=0\;.
\end{align*}
Since $v_\mathrm{left}$ is strictly positive, we conclude that $g$ is constant-0. That is, $g$ is an eigenfunction of $\Gamma$ for eigenvalue $\rho(\Gamma)$. By~\ref{en:EigenvaluesAndSimple}, that means that $g$ is a multiple of $v_\mathrm{right}$.
\end{proof}

An easy combination of Gelfand's formula and the nonnegativity of the leading eigenfunction is the following fact.
\begin{lem}\label{lem:rhomonotone}
Suppose that $\Gamma_1$ and $\Gamma_2$ are two digraphons with $\Gamma_1\le \Gamma_2$ (pointwise). Then $\rho(\Gamma_1)\le \rho(\Gamma_2)$.
\end{lem}
\begin{proof}
There is nothing to prove if $\rho(\Gamma_1)=0$, so assume $\rho(\Gamma_1)>0$. Let $f$ be an eigenfunction of $\Gamma_1$ corresponding to eigenvalue $\rho(\Gamma_1)$. By Theorem~\ref{thm:Schaefer74}\ref{en:LeadingEigenNonnegative}, $f$ is nonnegative. For every $k\in \N$, we have 
\[
\left(\frac{\Gamma_1}{\rho(\Gamma_1)}\right)^k f=f\;.
\]
Since $\Gamma_2^k\ge \Gamma_1^k$, and since $\Gamma_2^k$ and $f$ are nonnegative, we have
\[
\left(\frac{\Gamma_2}{\rho(\Gamma_1)}\right)^k f \ge f\;,
\]
and so the operator norm of $\left(\frac{\Gamma_2}{\rho(\Gamma_1)}\right)^k$ is at least~1. By Proposition~\ref{pro:Gelfand}, the spectral radius of $\frac{\Gamma_2}{\rho(\Gamma_1)}$ is at least~1, as was needed.
\end{proof}

\section{Proof of Theorem~\ref{thm:decompositionIntoComponents}}\label{sec:ProofDecomposition}
\subsection{Reduction to graphons}
We reduce the problem of the decomposition of digraphons into strong component to the problem of the decomposition of graphons into connected components, a concept introduced by Janson in~\cite{ConnectednessGraphons}. We recall Janson's terminology and one of his main results. Note that if $W$ is a graphon on $\Omega$, then for every $x\in\Omega$, we have $\degIn_W(x)=\degOut_W(x)$, and we call this quantity the \emph{degree} of $x$, \index{$\deg_W(x)$}$\deg_W(x)$.
\begin{defi}\label{def:connectedgraphon}
Let $I$ be a finite or countable set which does not contain~$0$. For a graphon $W$ on a probability space $(\Omega,\mu)$ we say that a decomposition $\Omega=X_0\sqcup \bigsqcup_{i\in I}X_i$ is a \emph{decomposition into connected components of $W$ with isolated elements $X_0$} if the following hold.
\begin{enumerate}[label=(\roman*)]
    \item\label{en:ZeroDegree} For almost every $x\in X_0$ we have $\deg_W(x)=0$.
    \item\label{en:Ah} For every $i\in I$ we have $\mu(X_i)>0$. If $A\sqcup B=X_i$ is a decomposition of $X_i$ into two sets of positive measure we have $\int_{A\times B}W>0$.
    \item\label{en:conn3} For every $i\in I$ we have $\int_{X_i\times (\Omega\setminus X_i)}W=0$.
\end{enumerate}
The sets $X_i$ ($i\in I$) are called the \emph{components} of $W$. We say that $W$ is \emph{connected} if $X_0$ is a nullset and $|I|=1$.
\end{defi}

One of the main results of~\cite{ConnectednessGraphons} is that a decomposition into connected components exists and is unique modulo nullsets. We state this result for later reference.
\begin{thm}\label{thm:ComponentsGraphon}
    Suppose that $W$ is a graphon on a probability space $\Omega$. Then there exists a decomposition $\Omega=X_0\sqcup \bigsqcup_{i\in I}X_i$ into connected components and isolated elements $X_0$ of the graphon $W$. Further, this decomposition is unique in the same sense as in Theorem~\ref{thm:decompositionIntoComponents}.
\end{thm}

Recall that in a finite graph, each two vertices in one connected component can be connected by a path. The next easy lemma is a graphon counterpart to this.
\begin{lem}\label{lem:pathsincomponentBETTER}
    Suppose that $U$ is a graphon on $\Omega$ and $X$ is a connected component in it. Let $g:\Omega^2\rightarrow\{0,1\}$ be the indicator of reachability, that is, $g(x,y)=1$ if and only if $y\in\ReachOut_U(x)$.\footnote{Note that since graphons are symmetric, this is equivalent to $y\in\ReachIn_U(x)$.} Then $g$ is constant-1 almost everywhere on $X$.
\end{lem}
\begin{proof}
We introduce a function $h:X\rightarrow[0,\mu(X)]$ as follows: for $x\in X$, let $h(x):=\int_{y\in X} g(x,y)\diff y$. We need to prove that $h$ is constant-$\mu(X)$ almost everywhere on $X$. For a contradiction, suppose otherwise. We distinguish two cases.
\begin{itemize}
\item \emph{$h$ is not a constant function.}\\
Then there exists $c\in(0,\mu(X))$ such that the sets $A:=h^{-1}([0,c))$ and $B:=h^{-1}([c,\mu(X)])$ have positive measure. As $X$ is connected, we have $\int_{A\times B}W>0$. In particular, the set $A^*:=\{x\in A:\deg_W(x,B)>0\}$ has positive measure. Take $y\in A^*$ arbitrary. On the one hand, $h(x)<c$ since $x\in A$. On the other hand, $x$ is connected to a positive measure of $y$'s with $h(y)\ge c$, and so $h(x)\ge c$ (see Lemma~\ref{lem:calculusofreachability}), a contradiction.
\item \emph{$h$ is a constant-$c$ function, for some $c\in [0,\mu(X))$}.\\
There exists $\eps>0$ so that the set $Z:=\{x\in X:\deg_W(x)\ge \eps\}$ we have $\mu(Z)\ge \eps$. Reachable vertices from a given vertex are a superset of the neighborhood of that vertex, which is of measure at least $\eps$ for every $x\in Z$. Hence, $h(x)\ge \eps$ for every $x\in Z$. We conclude that $c>0$.
Take $x\in X$ such that $h(x)=c$. Let $A\subseteq X$ be the set of vertices reachable from $x$, $\mu(A)=c$. We have that the sets $A$ and $B:=X\setminus A$ have positive measures. As $X$ is connected, we have $\int_{A\times B}W>0$. In particular, the set $B^*:=\{y\in B:\deg_W(y,A)>0\}$ has positive measure. But obviously, all these elements from $B^*$ are reachable from $x$, a contradiction.
\end{itemize}
\end{proof}

\begin{lem}\label{lem:pathsincomponent}
    Suppose that $U$ is a graphon on $\Omega$ and $X$ is a connected component in it. Suppose that $A\subseteq X$ is a set of positive measure. For each $k\in \N$, define
    \begin{equation}\label{eq:defc}
    F_k:=\left\{x_0\in X\::\: \int_{(x_1,\ldots,x_k)\in X^{k-1}\times A}\prod_{i=1}^k U(x_i,x_{i-1})>0\right\}\;.
    \end{equation}
    Then $\bigcup_{k=1}^\infty F_k=_0 X$.
\end{lem}
\begin{proof}
This follows from Lemma~\ref{lem:pathsincomponentBETTER}.
\end{proof}

With these preparations, we can begin our decomposition of the digraphon $\Gamma$. Define a graphon $W$ by
\begin{equation}\label{eq:symm}
W(x,y):=
\underbrace{\left(\sum_{\ell=1}^\infty 2^{-\ell}\Gamma^\ell(x,y)\right)}_{\textsf{(T1)}}
\underbrace{\left(\sum_{\ell=1}^\infty 2^{-\ell}\Gamma^\ell(y,x)\right)}_{\textsf{(T2)}}
\;.    
\end{equation}
Indeed, this definition gives a function which is nonnegative, bounded by~1, and symmetric.

Let $\Omega=X_0\sqcup \bigsqcup_{i\in I}X_i$ be the decomposition of $W$ into connected components given by Theorem~\ref{thm:ComponentsGraphon}. The main step of the proof of Theorem~\ref{thm:decompositionIntoComponents} is to show that $X_0$ is a fragmented set in $\Gamma$ and all $X_i$ ($i\in I$) are strong components. The next lemma about restricting $\Gamma$ to one connnected component $X_i$ is useful.

\begin{lem}\label{lem:RestrictionDoesNotDecrease}
    Suppose that $i\in I$ is given. Define a digraphon $\Gamma_*$ by $\Gamma_*:=\Gamma\llbracket X_i\rrbracket$. Then for each $\ell\in\N$ and almost all pairs $(x,y)\in X_i$ we have $\Gamma^\ell(x,y)=(\Gamma_*)^\ell(x,y)$.
\end{lem}
\begin{proof}
In the proof, we shall use the following trivial step: if $f$ and $g$ are two nonnegative functions on the same measure space, then
\begin{equation}\label{eq:ImplyMe}
    \int f=0\quad \mbox{implies}\quad\int fg=0\;.
\end{equation}

    Since we have $\Gamma\ge \Gamma_*$ pointwise, it is clear that $\Gamma^\ell(x,y)\ge(\Gamma_*)^\ell(x,y)$. So, to finish the proof, we need to prove that
\begin{equation}\label{KatieMelua}
    \int_{(x_0,x_1,\ldots,x_\ell)\in\Upsilon}\prod_{j=1}^\ell\Gamma(x_{j-1},x_j) =0\;,
\end{equation}
where $\Upsilon\subseteq X_i\times\Omega^{\ell-1}\times X_i$ is the set of those $(\ell+1)$-tuples which have at least one coordinate outside of $X_i$. Given $s,k\in \N$ with $s+k\le \ell$, let $\Upsilon_{s,k}\subseteq \Upsilon$ be the set of those $(\ell+1)$-tuples $(x_0,x_1,\ldots,x_\ell)$ for which the first $s$ many coordinates lie in $X_i$, the next $k$ many coordinates lie in $\Omega\setminus X_i$, the $(s+k+1)$st coordinate lies in $X_i$, and the remaining coordinates are unrestricted. Obviously, the sets $\{\Upsilon_{s,k}\}_{s,k\in \N, s+k\le \ell}$ cover the set $\Upsilon$. To prove~\eqref{KatieMelua}, it is hence enough to prove that for each $s,k\in \N$ with $s+k\le \ell$, $\int_{(x_0,x_1,\ldots,x_\ell)\in\Upsilon_{s,k}}\prod_{j=1}^\ell\Gamma(x_{j-1},x_j)=0$. To this end, we follow the proof scheme from~\eqref{eq:ImplyMe}. In particular, the function $f$ is set to $f:=\prod_{j=s}^{k+s}\Gamma(x_{j-1},x_j)$, while the function $g$ includes the remaining terms. Hence, the proof is finished by the following claim.
    \begin{claim}\label{claimA}
    For every $k\in\N$, we have $\int_{(y_0,y_1,\ldots,y_k)\in X_i\times (\Omega\setminus X_i)^{k-1}\times X_i}\prod_{j=1}^k \Gamma(y_{j-1},y_{j})=0$.
    \end{claim}
    \begin{proof}[Proof of Claim~\ref{claimA}]
    Suppose that this is not true. Let 
    \[
    Y:=\left\{y_1\in \Omega\setminus X_i\::\:\int_{(y_0,y_2,y_3,\ldots,y_k)\in X_i\times (\Omega\setminus X_i)^{k-2}\times X_i}\prod_{j=1}^k \Gamma(y_{j-1},y_{j})>0\right\}\;.
    \]
    Put equivalently, $Y$ consists of those $y_1\in \Omega\setminus X_i$ for which
    \begin{align}
\label{eq:Y1} \int_{y_0\in X_i}\Gamma(y_0,y_1)>0\quad\mbox{and}\\
\label{eq:Y2} \int_{(y_2,y_3,\ldots,y_k)\in (\Omega\setminus X_i)^{k-2}\times X_i}\prod_{j=2}^k \Gamma(y_{j-1},y_{j})>0\quad\mbox{.}
    \end{align}
    By the above assumption, $Y$ has positive measure. We claim that for every $y_1\in Y$ we have 
    \begin{equation}\label{eq:Wpos}
        \int_{x\in X_i}W(x,y_1)>0\;.
    \end{equation}
    Once we show this, we will get a contradiction to the fact that $X_i$ is a connected component. Indeed, $\int_{X_i\times(\Omega\setminus X_i)}W\ge\int_{X_i\times Y}W>0$, a contradiction with Definition~\ref{def:component}\ref{en:conn3}.

    To show~\eqref{eq:Wpos}, we will show for the given $y_1\in Y$, in the formula~\eqref{eq:symm} for $W(x,y_1)$,
    \begin{align}
\label{T1pos}
&\mbox{the term \textsf{(T1)} is positive for a positive measure of elements $x\in X_i$, and}\\
\label{T2pos}
&\mbox{the term \textsf{(T2)} is positive for almost all elements $x\in X_i$.}
\end{align}
Observe that~\eqref{eq:Y1} readily implies~\eqref{T1pos}. So, it remains to argue~\eqref{T2pos}. Let $A$ consist of those $y_k\in X_i$ for which 
$$\int_{(y_2,y_3,\ldots,y_{k-1})\in (\Omega\setminus X_i)^{k-2}}\prod_{j=2}^k \Gamma(y_{j-1},y_{j})>0\;.$$
By~\eqref{eq:Y2} has positive measure. Let $\{F_h\}_{h=1}^\infty$ be given by Lemma~\ref{lem:pathsincomponent} for the graphon $W$, the component $X_i$ and the set $A$. Now, observe that if $x\in X_i$ is such that it is contained in $F_h$, then $\Gamma^{(k-1)+h}(y_1,x)>0$. The way to think about this is that there is a positive density of paths of length $k-1$ from $y_1$ to each point of $A$, and from points from $A$ there is (aggregately, but not necessarily from each point individually) positive density of paths of length $h$ to $x$. Since $\bigcup_h F_h=_0 X$, we have that~\eqref{T2pos} follows.
\end{proof}
\end{proof}

\subsection{Sets $X_i$ are strong components}
Take $i\in I$. First, we shall show that $X_i$ is strongly connected in $\Gamma$. Let $A\sqcup B= X_i$ be as in Definition~\ref{def:component}\ref{en:StrongConn}. Suppose for contradiction that
\begin{equation}\label{eq:contradictioGamma}
    \int_{A\times B} \Gamma=0\;.
\end{equation}
Let $\{F_k\}_{k=1}^\infty$ be defined as in~\eqref{eq:defc} (for the graphon $W$, the component $X_i$ and the initial set $A$). Since $\mu(\cup_k F_k\setminus A)>0$, by sigma-additivity there exists $k\in \N$ so that $\mu(F_k\cap B)>0$. Extracting the term~\textsf{(T1)} from~\eqref{eq:symm}, we have that 
\[
\int_{(x_0,x_1,\ldots,x_k)\in A\times (X_i)^{k-1}\times (F_k\cap B)}\prod_{j=1}^k\left(\sum_{\ell_j=1}^\infty 2^{-\ell_j}\Gamma^{\ell_j}(x_{j-1},x_j)\right)>0\;.
\]
By Lemma~\ref{lem:RestrictionDoesNotDecrease}, we have that for the restriction $\Gamma_*$ of $\Gamma$ on the component $X_i$ that
\[
\int_{(x_0,x_1,\ldots,x_k)\in A\times (X_i)^{k-1}\times (F_k\cap B)}\prod_{j=1}^k\left(\sum_{\ell_j=1}^\infty 2^{-\ell_j}(\Gamma_*)^{\ell_j}(x_{j-1},x_j)\right)>0\;.
\]
By sigma-additivity of the integral, there exist numbers $\ell_1,\ldots,\ell_k\in\N$ so that
\begin{equation}\label{AmyMcDonald:ThisIsTheLife}
\int_{(x_0,x_1,\ldots,x_k)\in A\times (X_i)^{k-1}\times (F_k\cap B)}\prod_{j=1}^k (\Gamma_*)^{\ell_j}(x_{j-1},x_j)>0\;.    
\end{equation}
We now expand the terms $(\Gamma_*)^{\ell_j}$; each such term involves an internal $(\ell_j-1)$-dimensional integration over $X_i$ because Definition~\ref{def:digraphonpower} involves~\eqref{eq:densPathTwoTerminals}. Hence, we can rewrite~\eqref{AmyMcDonald:ThisIsTheLife} as
\begin{equation}\label{AmyMcDonaldFinished}
\int_{(x_0,x_1,\ldots,x_{\sum_j \ell_j})\in A\times (X_i)^{\sum_j \ell_j-1}\times (F_k\cap B)}\prod_{t=1}^{\sum_j \ell_j} \Gamma_*(x_{t-1},x_t)>0\;.    
\end{equation}
For $L\in\N$, let $\Upsilon_L$ consist of those tuples $(x_0,x_1,\ldots,x_{\sum_j \ell_j})\in A\times (X_i)^{\sum_j \ell_j-1}\times (F_k\cap B)$ for which the first $L+1$ many components lie in $A$ and the $(L+1)$-st lies in $F_k\cap B$. Obviously, the sets $\{\Upsilon_L\}_{L=1}^\infty$ cover $A\times (X_i)^{\sum_j \ell_j-1}\times (F_k\cap B)$. For a given $L\in\N$, \eqref{eq:contradictioGamma} can be rewritten as 
\begin{equation*}
\int_{(x_0,x_1,\ldots,x_{\sum_j \ell_j})\in \Upsilon_L}\Gamma_*(x_{L},x_{L+1})=0\;.    
\end{equation*}
Summing over all $L$ and using the proof scheme~\eqref{eq:ImplyMe} gives a contradiction with~\eqref{AmyMcDonaldFinished}. So, $X_i$ is indeed strongly connected.

Let us now verify Definition~\ref{def:component}\ref{en:StrongComp}. Suppose that $i\in I$ and $Y\subseteq \Omega$ is such that $\mu(X_i\cap Y),\mu(Y\setminus X_i)>0$. We need to show that $Y$ is not strongly connected. Let $A:=X_i\cap Y$. Let $B^-:=\ReachOut_\Gamma(A)\cap (Y\setminus X_i)$, and $B^+:=\ReachIn_\Gamma(A)\cap (Y\setminus X_i)$.
\begin{claim}\label{claim:BNull}
The set $B^-\cap B^+$ is null. 
\end{claim}
\begin{proof}[Proof of Claim~\ref{claim:BNull}]
Indeed, it is obvious from~\eqref{eq:symm} that for each $b\in B^-\cap B^+$ we have $\deg_W(b,A)>0$. Recall that $B^-\cap B^+$ is disjoint from the component $X_i$, and thus $\int_{(B^-\cap B^+)\times X_i} W=0$. The claim follows.
\end{proof}

To verify Definition~\ref{def:component}\ref{en:StrongComp}, we need to find a partition $Y=A^*\sqcup B^*$ into sets of positive measures so that $\int_{A^*\times B^*}\Gamma=0$. In order to ensure that $A^*$ and $B^*$ have positive measures, we need to distinguish three cases.

\begin{itemize}
    \item If $\mu(B^+)>0$, then set $A^*:=Y\setminus B^+$, $B^*:=B^+$. The set $A^*$ contains $A$, and so $A^*$ and $B^*$ have positive measures. Suppose for contradiction that $\int_{A^*\times B^*}\Gamma>0$. That means that there is a set $Z\subseteq A^*$ of positive measure so that for every $z\in Z$ we have $\degOut_\Gamma(z,B^+)>0$. But each such $z$ lies in $\ReachIn_\Gamma(A)$, a contradiction.
    \item If $\mu(B^-)>0$, then set $A^*:=B^-$, $B^*:=Y\setminus B^-$. The argument is the same as in the previous case.
    \item If $\mu(B^+)=\mu(B^-)=0$ then set $A^*:=A$ and $B^*:=Y\setminus A$. $A^*$ and $B^*$ have positive measures. It is clear that $\int_{A^*\times B^*}\Gamma=0$.
\end{itemize}

\subsection{Set $X_0$ is fragmented}
Let $Y\subseteq X_0$ of positive measure be as in Definition~\ref{def:component}\ref{en:Fragmented}. For $x\in Y$, let $f(x)\in[0,\mu(Y)]$ be the measure of the set of points in $Y$ which can be reached from $x$. We distinguish three cases.
\begin{itemize}
    \item \emph{$f$ is constant-0 almost everywhere on $Y$.}\\
    In that case we partition $Y=A\sqcup B$ into two arbitrary sets of positive measure. Obviously, we have $\int_{A\times B}\Gamma=0$. Definition~\ref{def:component}\ref{en:Fragmented} is satisfied.
    \item \emph{$f$ is constant-$c$ almost everywhere on $Y$, for some $c>0$.}\\
    We will show a contradiction. Let us take $x\in Y$ with $f(x)=c$. Let $Z\subseteq Y$ be those points that can be reached from $x$, $\mu(Z)=c$. We claim that for almost every point $z\in Z$, the set of points that can be reached from $z$ is equal to $Z$ up to a nullset. Indeed, suppose otherwise. That is, there is a set $C\subseteq Z$ of positive measure of points $z$ with $\mu(\ReachOut_\Gamma(z)\setminus Z)>0$. But the Lemma~\ref{lem:calculusofreachability} tells us that $\mu(\ReachOut_\Gamma(x)\setminus Z)>0$, a contradiction.
    
    We concluded that for almost all the pairs $(z,z')\in Z^2$, $z$ can be reached from $z'$ and $z'$ can be reached from $z$. In particular $W(z,z')>0$. So, the degrees of vertices in $Z$ in the graphon $W$ are positive, a contradiction with Definition~\ref{def:connectedgraphon}\ref{en:ZeroDegree}.
    
    \item \emph{$f$ is not constant almost everywhere on $Y$.}\\
    Take $c\in [0,\mu(Y))$ such that $A:=f^{-1}([0,c])$ and $B:=f^{-1}((c,\mu(Y)])$ are two sets of positive measure. We claim that $\int_{A\times B}\Gamma=0$. Indeed, by Lemma~\ref{lem:calculusofreachability}, for every $x\in A$, we have $\int_{y\in B}\Gamma(x,y)=0$. Hence, Definition~\ref{def:component}\ref{en:Fragmented} is satisfied.
\end{itemize}

\subsection{Uniqueness}
Suppose that $\{X'_i\}_{i\in I'\cup\{0\}}$ is another decomposition of $\Gamma$ into strong components and a fragmented set. First, observe that Definition~\ref{def:component}\ref{en:StrongComp} implies that if for some $i\in I$ and $i'\in I'$ we have that $\mu(X_i\cap X'_{i'})>0$ then $X_i=_0 X'_{i'}$. That is, the only thing left to rule out are potential strong components $X_i$ with the property $X_i\subseteq X'_0$ or $X'_{i'}$ with the property $X'_{i'}\subseteq X_0$. But these cannot exist by the antagonistic properties of Definition~\ref{def:component}\ref{en:StrongComp} and Definition~\ref{def:component}\ref{en:Fragmented}.

\subsection{Useful additional properties}\label{ssec:AdditionalDigraphDecomposition}
Given the fundamental nature of the decomposition into strong components and the fragmented set, here we establish additional useful properties of this decomposition. Proposition~\ref{prop:reachabilityFromAnywhereToAnywhere} asserts that almost all pairs of points in a strong component are connected in both directions. In Proposition~\ref{prop:CondensationDigraph} we propose a counterpart of the well-known notion of the `condensation digraph' for digraphons. In subsequent comments, we extend this to `enhanced condensation digraph' which captures also connectivity properties outside of the fragmented set. Proposition~\ref{prop:X0main} deals with the fragmented set. The original Definition~\ref{def:component}\ref{en:Fragmented} is not a very descriptive one. In Proposition~\ref{prop:X0main} we for example show that the fragmented set can be linearly ordered so that for $x$ and $y$ in the fragmented set, we have $\Gamma(x,y)>0$ only if $x$ precedes $y$. We also show that the spectral radius of the fragmented set is~0. Last, in Proposition~\ref{prop:cyclesconfined}, we show that all directed cycles are confined to strong components.

\begin{prop}\label{prop:reachabilityFromAnywhereToAnywhere}
Suppose that $X_i$ is a strong component of a digraphon $\Gamma$. Then for almost every $x\in X_i$, almost all points of $X_i$ are reachable from $x$, and for almost every $x\in X_i$, $x$ can be reached from almost every point of $X_i$.
\end{prop}
\begin{proof}
We proved that strong components of $\Gamma$ correspond to components of $W$. Also, from the way $W$ was defined in~\eqref{eq:symm}, we have that if $x$ and $y$ are reachable in $W$ then $x$ is reachable from $y$ in $\Gamma$ and also that $y$ is reachable from $x$ in $\Gamma$. Hence, the claim follows from Lemma~\ref{lem:pathsincomponentBETTER}.
\end{proof}

Recall that decomposing a finite digraph into its strong components leads to the notion of a condensation digraph. This is a digraph whose vertices correspond to the strong components, and in which a directed edge between two such vertices indicates the existence of at least one directed edge between the corresponding components in the original digraph. A crucial property of a condensation digraph is that it is acyclic. Here, we establish an analogue for digraphons.
\begin{prop}\label{prop:CondensationDigraph}
Suppose that $\Gamma$ is a digraphon and $\Omega=X_0\sqcup \bigsqcup_{i\in I} X_i$ is its decomposition into strong components. Create a (finite or countable) digraph on $D$ vertex set $I$ as follows. A pair of distinct vertices $(i,j)\in I^2$ is a directed edge if and only $\Gamma$ is not constant-0 on $X_i\times X_j$. We call $D$ the \emph{condensation digraph of $\Gamma$}.
\begin{enumerate}[label=(\roman*)]
\item\label{en:CondensationReachability} Suppose that $\ell\in \N$ and $i_1i_2\ldots i_\ell$ is a directed walk in $D$. Then for almost every pair of vertices $(x,y)\in X_{i_1}\times X_{i_\ell}$, we have that $y$ is reachable from $x$ in $\Gamma$.
\item\label{en:CondensationAcyclic} The condensation digraph $D$ is acyclic, that is, it contains no finite directed cycle. In particular, applying this for directed cycles of length~2, the digraphon $\Gamma$ is constant-0 on at least one of $X_i\times X_j$ and $X_j\times X_i$ for each pair of distinct $i,j\in I$.
\end{enumerate}  
\end{prop}
\begin{proof}
We first prove~\ref{en:CondensationReachability}.
We prove the claim by induction on $\ell$. First, we deal with the base case $\ell=1$. That is, we need to prove that for almost every pair $(x,y)\in X_{i_1}\times X_{i_1}$, $x$ is reachable from $y$. This is the subject of Proposition~\ref{prop:reachabilityFromAnywhereToAnywhere}.

Let us now deal with the induction step. Let $i_1i_2\ldots i_{\ell+1}$ be a directed walk in $D$. We know that $\Gamma$ is not constant-0 on $X_{i_\ell}\times X_{i_{\ell+1}}$. In particular, there is a set $Z\subseteq X_{i_{\ell+1}}$ of positive measure, each element of which has positive indegree from the set $X_{i_{\ell}}$. By the induction hypothesis, for almost all pairs $(x,y)\in X_{i_1}\times X_{i_{\ell}}$ we have that $y$ is reachable from $x$. By Lemma~\ref{lem:calculusofreachability}, for almost all pairs $(x,z)\in X_{i_1}\times Z$, $z$ is reachable from $x$. Combined with Proposition~\ref{prop:reachabilityFromAnywhereToAnywhere}, we conclude that for almost all $(x,w)\in X_{i_1}\times X_{i_{\ell+1}}$, $w$ is reachable from $x$, as was needed.

We can now prove~\ref{en:CondensationAcyclic}.
Suppose that $C$ is a directed cycle in the condensation digraph. Then we can take two vertices of it, say, $i$ and $j$, and use $C$ to find a directed path from $i$ to $j$ and from $j$ to $i$. Part~\ref{en:CondensationReachability} tells us that for almost every pair $(x,y)\in X_i\times X_j$ we have that $x$ is reachable from $y$ and $y$ is reachable from $x$. That means that in~\eqref{eq:symm}, we have $W(x,y)>0$. But then $X_i$ and $X_j$ are not different connected components of $W$, a contradiction.
\end{proof}

\begin{figure}\centering
	\includegraphics[scale=0.7]{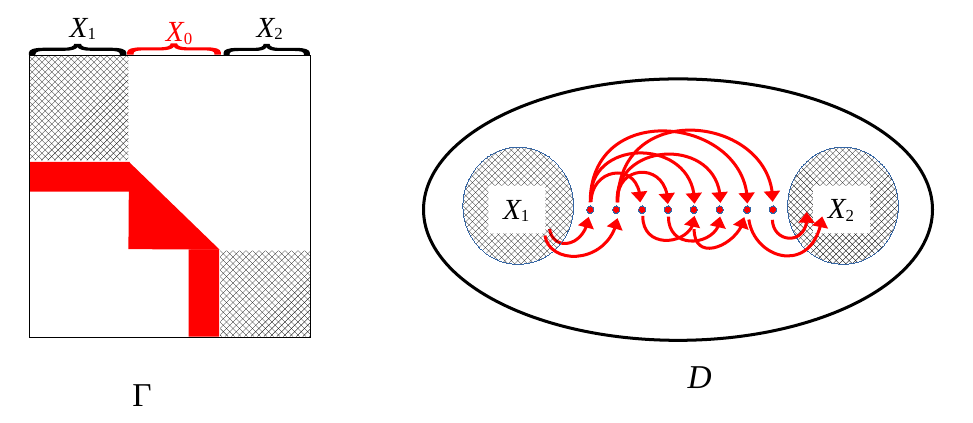}
	\caption{An example of a digraphon $\Gamma$, and digraph $D$ which approximates $\Gamma$. The condensation digraph of $\Gamma$ consists of two isolated vertices. Also, note that any directed path from $X_1$ to $X_2$ needs to go through several vertices in $X_0$.}
	\label{fig:ConnectivitzViaX0}
\end{figure}

Proposition~\ref{prop:CondensationDigraph} introduced the bare condensation digraph on the strong components $\{X_i\}_{i\in I}$ but did not treat the set $X_0$. For finite digraphs, Proposition~\ref{prop:CondensationDigraph}\ref{en:CondensationReachability} can be reversed: if $u\in C_i$ and $v\in C_j$ are vertices in two strongly connected components $C_i$ and $C_j$ of a digraph, and there is no directed walk from $C_i$ to $C_j$. This is not true for the condensation digraph. Indeed, Figure~\ref{fig:ConnectivitzViaX0} shows a digraphon with two strong components $X_1$ and $X_2$ whose condensation digraph consists of two isolated vertices. Points of $X_2$ can be reached from $X_1$, though only using the fragmented sets $X_0$. Hence, for some scenarios, there is a more telling notion of \emph{extended condensation digraph} $D^+$. $D^+$ is a supergraph of $D$ on the same vertex set $I$. A pair $(i,j)$ which is not an edge is added to $D^+$ if there is a positive density of paths (of any length) from $X_i$ to $X_j$ whose internal vertices are only in $X_0$. It is easy to show, that with this definition, Proposition~\ref{prop:CondensationDigraph}\ref{en:CondensationReachability} can be reversed: if $(i,j)\in I^2$ are such that there is a directed walk from $i$ to $j$ in $D^+$, then for almost every pair of vertices $(x,y)\in X_i\times X_j$, $y$ is reachable from $x$. If, on the other hand, such a directed walk from $i$ to $j$ does not exist, then the assertion holds for almost no pair $(x,y)\in X_i\times X_j$.

Proposition~\ref{prop:X0main} is our main technical statement about the fragmented set in a decomposition of a digraphon into strong components. While the first three items may seem little technical, the last two show how they are useful. In particular, Item~\ref{en:X0linearorder} tells us that the set $X_0$ (or an arbitrary subset of it) can be linearly ordered so that there are no edges going from bigger elements to smaller elements.

\begin{prop}\label{prop:X0main}
Suppose that $\Gamma$ is a digraphon and $\Omega=X_0\sqcup \bigsqcup_{i\in I} X_i$ is its decomposition into strong components. Let $X\subseteq X_0$ be arbitrary of positive measure. Let $g:X\to [0,\mu(X)]$ be defined by $g(x):=\mu(\ReachOut_\Gamma(x)\cap X)$. For $a\in[0,1]$, let $F_a:=g^{-1}([a,1])$.
\begin{enumerate}[label=(\roman*)]
\item\label{en:X0main1} Let $a\ge 0$ be arbitrary. Then for almost every $x\in F_{a}$ we have $\mu(\ReachOut_\Gamma(x)\cap X\setminus F_{a})\ge a$.
\item\label{en:X0main2} Let $a\ge 0$ be arbitrary. Then we have $\mu(F_a)\le \mu(X)-a$.
\item\label{en:X0linearorder} $\Gamma$ is constant-0 almost everywhere on the set $\{(x,y)\in X^2: g(x)\ge g(y)\}$.
\item\label{en:X0SpectralRadius} We have $\rho(\Gamma\llbracket X\rrbracket)=0$.
\end{enumerate}
\end{prop}

For the proof, we shall need a lemma which we might call `a qualitative Caccetta--Häggkvist Conjecture for digraphons'. Recall that the Caccetta--Häggkvist Conjecture assert that if $D$ is an $n$-vertex digraph with minimum outdegree at least $n/k$ for some $k\in \N$, then $D$ contains an oriented cycle of length at most $k$. While the question is wide open, qualitative relaxations are known. As an example which will be useful for us, Shen~\cite{Shen2002} proved that any such digraph contains a cycle of length at most $k+73$. Below we give a similar result for digraphons. We do not require any bound on the length of the cycles (even though our proof easily gives one), but we require positive density of such cycles (which is the only reasonable counterpart to the existence of a subdigraph in finite digraphs).
\begin{lem}\label{lem:CaccettaHaggkvist}
Suppose that $\Psi$ is a digraphon on ground space $\Xi$ with $\degminOut(\Psi)>0$. Then there exists $k\ge 2$ so that $t(C_k,\Psi)>0$.
\end{lem}
\begin{proof}
This is almost exactly the result of Shen mentioned above. The only difference is that Shen's result is for finite digraphs whereas our claim is about digraphons. Our derivation is standard, using the Regularity lemma.

Suppose that the probability measure on $\Xi$ is $\xi$. Let $a:=\degminOut(\Psi)$. We apply Theorem~\ref{thm:Regularity}, with error parameter $\eps=(d/4)^{2000\cdot d^{-2}}$ and density $d:=a/3$ on the digraphon $\Psi$. Theorem~\ref{thm:Regularity} outputs a partition $\Xi=Z_1\sqcup \cdots \sqcup Z_t$ and a digraph $D$ on vertex set $[t]$. The minimum outdegree of $D$ is at least $at/2$.
For $i,j\in [t]$ we will refer to densities $p_{i,j}$ and digraphons $\Psi^{avg}_{i,j}$ from Definition~\ref{def:clusterdigraph}. Shen's result~\cite{Shen2002} asserts that $D$ contains a directed cycle $C=v_1v_2\ldots v_k$ of length $k\le \frac{2}{d}+73$. The product $p:=\prod_{i=1}^k p_{v_i,v_{i+1}}$ (with the notation $v_{k+1}=v_1$) satisfies 
\begin{equation}\label{eq:lbp}
p\ge (d/3)^{\frac{2}{d}+73}\;.
\end{equation}

We recall Definition~\ref{def:restriction}. Let $Z:=\cup_{\ell=1}^k Z_{v_\ell}$, and $\Gamma:=\Psi_{\restriction Z\times Z}$, and $\Gamma^{avg}:=\left(\sum_{\ell=1}^k \Psi^{avg}_{v_\ell,v_{\ell+1}}\right)_{\restriction Z\times Z}$. That is, $\Gamma$ zooms on the part of $\Psi$ containing cycle $C$, and $\Gamma^{avg}$ is the averaged version. By summing up $k$ many times the cut distances in Definition~\ref{def:clusterdigraph}\ref{en:clCut} (and taking into account the rescaling of the measure by a factor of $\frac{t}{k}$ from Definition~\ref{def:restriction}), we have $\cutnormdist(\Gamma,\Gamma^{avg})<k\cdot (\eps\cdot t^{-2})\cdot(\frac{t}{k})^2=\eps/k$.

The defining product for $p$ is what appears in~\eqref{eq:defHomDens} for $t(C_k,\Gamma^{avg})$ except that in the latter only when the elements tuple $(x_1,\ldots,x_k)$ are confined to consecutive cells of the partition $Z=Z_{v_1}\sqcup Z_{v_2}\sqcup \ldots\sqcup Z_{v_\ell}$. Hence, 
\[
t(C_k,\Gamma^{avg})\ge p\cdot (1/k)^k
\geByRef{eq:lbp}
(d/3)^{\frac{2}{d}+73}\cdot (1/k)^k
>\eps\;.
\]
Hence, by Lemma~\ref{lem:CountingLemma}, $t(C_k,\Gamma)\ge t(C_k,\Gamma^{avg})-e(C_k)\cdot \cutnormdist(\Gamma,\Gamma^{avg})>0$. Since $\Gamma$ is contained in a rescaled version of $\Psi$, we conclude that $t(C_k,\Psi)>0$.
\end{proof}

\begin{proof}[Proof of Proposition~\ref{prop:X0main}\ref{en:X0main1}]
Suppose that this fails for some $a$. That means that there exists $\eps>0$ so that the set
\[
Z:=\left\{x\in F_a\::\: \mu(\ReachOut_\Gamma(x)\cap X\setminus F_{a})\le a-\eps\right\}
\]
has positive measure. By Lemma~\ref{lem:calculusofreachability}, almost every vertex $y$ reachable from a given vertex $x\in Z$ has itself $\mu(\ReachOut_\Gamma(y)\cap X\setminus F_{a})\le a-\eps$. In particular, if $y$ is in addition in $F_a$, then $y\in Z$. Since $x\in Z\subseteq F_a$, we know that 
\begin{equation}\label{eq:ReachMe}
\mu(\ReachOut_\Gamma(x)\cap Z)\ge\eps
\;.
\end{equation}

Construct a digraphon $\Psi$ on ground space $Z$ equipped with the natural rescaling of measure $\mu$, $\mu_Z(\cdot):=\frac{\mu(\cdot\cap Z)}{\mu(Z)}$, as follows. For $(x,y)\in Z^2$, $\Psi(x,y)$ is~1 or~0, depending on whether $y\in\ReachOut_\Gamma(x)$ or not. Then, by~\eqref{eq:ReachMe}, we have $\degminOut(\Psi)\ge \frac{\eps}{\mu(Z)}$ (the denominator takes care of rescaling between $\mu$ and $\mu_Z$). By Lemma~\ref{lem:CaccettaHaggkvist}, we have $t(C_k,\Psi)>0$ for some $k\ge 2$. That is, the set $P\subseteq Z^2$ of pairs $(x_1,x_2)\in Z^2$ such that, writing $x_{k+1}=x_1$,
\[
\int_{x_3\in Z}\int_{x_4\in Z}\ldots\int_{x_k\in Z}\prod_{i=1}^k\Psi(x_i,x_{i+1})>0
\]
has positive measure (either with respect to $\mu$ or with respect to $\mu_Z$, this is equivalent since we only care about positivity). For each such $(x_1,x_2)\in P$ we have that $x_2\in \ReachOut_\Gamma(x_1)$ by the definition of $\Psi$. But also, using paths via $x_3x_4\ldots x_k$, we see that $x_1\in \ReachOut_\Gamma(x_2)$. This means that in~\eqref{eq:symm}, we have $W(x_1,x_2)>0$. As $P^2\subseteq X^2$, we deduced that $W$ is not constant-0 almost everywhere on $X^2$, a contradiction.
\end{proof}
\begin{proof}[Proof of Proposition~\ref{prop:X0main}\ref{en:X0main2}]
Set $G_a:=X\setminus F_a$. The proof proceeds by double-counting pairs $(x,y)\in F_a\times G_a$ such that (*1)~$y\in \ReachOut_\Gamma(x)$ or equivalently, (*2)~$x\in \ReachIn_\Gamma(y)$. To get a lower-bound on the measure of $y$'s that satisfy~(*1) for a given $x$, we use Part~\ref{en:X0main1}. To get an upper-bound on the measure of $x$'s that satisfy~(*2) for a given $y$, we use that $\mu(\ReachIn_\Gamma(y)\cap F_a)\le \mu(F_a)$. That is,
\[
\mu(F_a)a
\le 
\mu^2
\left(\left\{(x,y)\in F_a\times G_a\::\:\mbox{$(x,y)$ satisfies (*1) and/or (*2)}\right\}\right)
\le 
\mu(F_a)\mu(G_a)\;.
\]
This yields $\mu(F_a)=0$ or $\mu(G_a)\ge a$. Either case suffices to conclude the statement.
\end{proof}

\begin{proof}[Proof of Proposition~\ref{prop:X0main}\ref{en:X0linearorder}]
We need to check that for every $c\in [0,\mu(X)]$, $\Gamma$ is constant-0 on $g^{-1}([c,\mu(X)])\times g^{-1}([0,c])$. This follows by Lemma~\ref{lem:calculusofreachability} for the set $g^{-1}([c,\mu(X)))\times g^{-1}([0,c])\cup g^{-1}([c,\mu(X)])\times g^{-1}([0,c))$. So, in the rest we prove that $\Gamma$ is constant-0 on $Z\times Z$, where $Z:=g^{-1}(\{c\})$. Assume that $Z$ has positive measure as there is nothing to prove otherwise.
For $x\in Z$, we decompose the value $g(x)$, $g(x)=f(x)+h(x)$, $f(x):=\mu(\ReachOut_\Gamma(x)\cap Z)$, $h(x):=\mu(\ReachOut_\Gamma(x)\cap X\setminus Z)$. If $f$ is constant-0 on $Z$, then $\Gamma$ is constant-0 on $Z\times Z$ and we are done. If $f$ is constant-$c'$ for some $c'>0$, we get a contradiction with Lemma~\ref{lem:CaccettaHaggkvist} in the same way as we did in the proof of Part~\ref{en:X0main1}, this time on the ground set $Z$. 
Hence, it remains to deal with the case that for some $c''>0$ we have the sets $S_1:=f^{-1}([0,c''])$ and $S_2:=f^{-1}((c'',c])$ have positive measures. Note that equivalently, $S_1:=h^{-1}([c-c'',c])$ and $S_2:=h^{-1}([0,c-c''))$. Lemma~\ref{lem:calculusofreachability} tells us that for $x\in S_2$, $\ReachOut_\Gamma(x)\cap S_1$ is null. Hence,
\[
\mu(\ReachOut_\Gamma(x)\cap S_2)=\mu(\ReachOut_\Gamma(x)\cap Z)=f(x)>c''\;.
\]
Therefore, we get a contradiction with Lemma~\ref{lem:CaccettaHaggkvist} in the same way as we did in the proof of Part~\ref{en:X0main1}, this time on the ground set $S_2$. 
\end{proof}

\begin{proof}[Proof of Proposition~\ref{prop:X0main}\ref{en:X0SpectralRadius}]
Rescale $X$ to be a probability space, $\mu_X(\cdot):=\frac{\mu(\cdot)}{\mu(X)}$. Consider the lower-triangular digraphon $T\in L^\infty (X^2)$, defined by $T(x,y)=\mathbbm{1}_{\{g(x)< g(y)\}}$. Let $B:X\to [0,\mu(X)]$ be defined by $B(t):=\mu_X(F_t)$. Simple inductive calculations show that for each $k\in \N$ and for each $x,y\in[0,1]$ we have
\[
T^k(x,y)\le \frac{1}{(k-1)!}\cdot \mathbbm{1}_{\{g(x)< g(y)\}} (B(y)-B(x))^k\;.
\]
Since, in general $\|\cdot \|_{\mathrm{op}}\le \|\cdot\|_\infty$, we have $\|T^k \|_{\mathrm{op}}\le \frac{1}{(k-1)!}$ for each $k\in\N$. By Proposition~\ref{pro:Gelfand}, we have $\rho(T)\le \lim_{k\to\infty}\sqrt[k]{\frac{1}{(k-1)!}}=0$. By the monotonicity of the spectral radius (Lemma~\ref{lem:rhomonotone}), we have $\rho\left(\Gamma_{\restriction X\times X}\right)\le \rho(T)=0$. Hence, $\rho\left(\Gamma\llbracket X\rrbracket\right)=0$.
\end{proof}

Last useful property suggests that directed cycles are confined to strong components.
\begin{prop}\label{prop:cyclesconfined}
Suppose that $\Gamma$ is a digraphon and $\Omega=X_0\sqcup \bigsqcup_{i\in I} X_i$ is its decomposition into strong components. Let $k\in \{2,3,\ldots\}$ be arbitrary. Then we have $t(C_k,\Gamma)=\sum_{i\in I}t(C_k,\Gamma\llbracket X_i\rrbracket)$.
\end{prop}
\begin{proof}
We need to prove that in the integral $\int_{x_1,\ldots,x_k}$ defining $t(C_k,\Gamma)$ there is no contributions of tuples $(x_1,x_2,\ldots x_k)\in X_{i_1}\times X_{i_2}\times \ldots \times X_{i_k}$ with $\{i_1,i_2,\ldots,i_k\}\subseteq I$ being not all the same, and also that there is no the contribution of tuples $\bigcup_{\ell=1}^k\{(x_1,x_2,\ldots x_k):x_\ell\in X_0\}$.

As for the first part, we see that this is clearly the case for the given tuple $(i_1,i_2,\ldots,i_k)$ as above if for some $s\in [k]$ we have that $i_si_{s+1}$ (using the cyclic notation) is not a directed edge of the condensation digraph of our partition. But Proposition~\ref{prop:CondensationDigraph} tells us that such a nonedge always exists.

As for the second part, suppose for contradiction that $Z=\{x_1\in X_0:t^{\bullet}_{x_1}(C_k^{\bullet},\Gamma)>0\}$ has positive measure. Let $Y\subseteq \Omega$ be defined as those points $x_2$ where 
\[
\int_{x_1\in Z,x_3\in \Omega,x_4\in \Omega,\ldots, x_k\in \Omega}
\prod_{s\in [k]}\Gamma(x_i,x_{i+1})>0\;.
\]
Obviously, $Y$ has positive measure. We see that the above definition of $Y$ is equivalent to
\[
Y=\{
x_2\in \Omega:\mu(\{x_1\in Z:\Gamma(x_1,x_2)>0\mbox{ and }\Gamma^{k-1}(x_2,x_1)\})>0)
\}\;.
\]
We now recall the definition of the graphon $W$ in~\eqref{eq:symm}. In particular, the above tells us that $\int_{Z\times Y} W>0$, a contradiction to the fact that $X_0$ was created as elements of zero degree in $W$.
\end{proof}

\section{Spectral radius of digraphons and strong components}\label{ssec:SpectralRadiusAndStrongComponents}
The main result of this section, Proposition~\ref{prop:spectralradiusAndStrongComponents}, connects spectral properties of a digraphon with the structure of its strong components.
Here, we work with the convention that the supremum of an empty set is~0 (rather than the more common $\sup\emptyset=-\infty)$.
\begin{prop}\label{prop:spectralradiusAndStrongComponents}
Suppose that $\Gamma$ is a digraphon on $\Omega$. Let $\Omega=X_0\sqcup \bigsqcup_{i\in I}X_i$ be its decomposition into strong components.
\begin{enumerate}[label=(\roman*)]
\item\label{en:spectralradiusMaximum} We have $\rho(\Gamma)=\sup_{i\in I}\rho(\Gamma\llbracket X_i\rrbracket)$.
The supremum $\sup_{i\in I}\rho(\llbracket X_i\rrbracket)$ is in fact a maximum. 
\item\label{en:CharSpectralRadius0} We have $\rho(\Gamma)=0$ if and only if $I=\emptyset$.
\end{enumerate}
\end{prop}
\begin{rem}\label{rem:spectralradiusAndStrongComponentsUnbounded}
As can be seen from the proof below, Proposition~\ref{prop:spectralradiusAndStrongComponents} holds even in a somewhat more general setting. Namely, suppose that $p\in[1,\infty)$, and that $\Gamma$ is nonnegative and such that as an integral kernel operator on $L^p(\Omega)$ it is bounded and compact. In that case, the correct key formula in Proposition~\ref{prop:spectralradiusAndStrongComponents} is $\rho_p(\Gamma)=\sup_{i\in I}\rho_p(\Gamma\llbracket X_i\rrbracket)$, where $\rho_p$ is the spectral radius when $\Gamma$ is viewed as an operator on $L^p(\Omega)$. This extension is useful in~\cite{HladkySavicky:Inhomo2SAT}.
\end{rem}
\begin{proof}[Proof of Proposition~\ref{prop:spectralradiusAndStrongComponents}\ref{en:spectralradiusMaximum}]
The case $I=\emptyset$ is handled by Proposition~\ref{prop:X0main}\ref{en:X0SpectralRadius}. We will assume that $I\neq \emptyset$.

The inequality $\rho(\Gamma)\ge \sup_{i\in I}\rho(\Gamma\llbracket X_i\rrbracket)$ is trivial. Indeed, each digraphon on the right-hand side is a sub-digraphon of $\Gamma$, and so the inequality follows from monotonicity of the spectral radius (Lemma~\ref{lem:rhomonotone}).
The proof of the inequality $\rho(\Gamma)\le \sup_{i\in I}\rho(\Gamma\llbracket X_i\rrbracket)$ has two steps. First, we prove the proposition when $I$ is finite. Then we use the finite case to cover the infinite case as well.

\underline{The case when $I$ is finite.} Consider the condensation digraph $D$ on $I$ as in Proposition~\ref{prop:CondensationDigraph}. Order the vertices $I$ as $I=\{t_1,\ldots,t_\ell\}$ so that for each $i\in [\ell]$ there are no directed edges going to $t_i$ from $\{t_{i+1},\ldots,t_\ell\}$. This is possible as $D$ is acyclic (Proposition~\ref{prop:CondensationDigraph}\ref{en:CondensationAcyclic}). Further, for each $i\in\{0,1,\ldots,\ell\}$, let $Z_i:=X_0\setminus \ReachOut_\Gamma\left(\bigcup_{j=i+1}^\ell X_{t_j}\right)$. We have $Z_0\subseteq Z_1\subseteq \ldots \subseteq Z_\ell=X_0$. 

The core of the proof is the following claim.
\begin{claim}\label{claimB}
Suppose that $f$ is an eigenfunction of $\Gamma$ corresponding to eigenvalue $\nu$. Then the following is true.
\begin{enumerate}[label=(\roman*)]
    \item\label{enCl:1} $f$ is constant-0 on $Z_0$.
    \item\label{enCl:2} If for some $i\in[\ell]$, $f$ is constant-0 on $Z_{i-1}\cup \bigcup_{j=1}^{i-1}X_{t_j}$, then either $f$ is constant-0 on $X_{t_i}$ or $\mathbbm{1}_{X_{t_i}}f$ is an eigenfunction of $\Gamma\llbracket X_{t_i}\rrbracket$ corresponding to eigenvalue $\nu$.
    \item\label{enCl:3} If some some $i\in[\ell]$, $f$ is constant-0 on $Z_{i-1}\cup \bigcup_{j=1}^{i}\cup X_{t_j}$, then $f$ is also constant-0 on $Z_i$.
\end{enumerate}
\end{claim}
\begin{proof}[Proof of Claim~\ref{claimB}]
All the parts of the claim follow the same idea, namely that the chain 
\[
\emptyset
\;,\;
Z_0 
\;,\;
Z_0\cup X_{t_1}
\;,\;
Z_1\cup X_{t_1}
\;,\;
Z_1\cup X_{t_1}\cup X_{t_2}
\;,\;
Z_2\cup X_{t_1}\cup X_{t_2}
\;,\; \ldots \;,\;
Z_\ell\cup \bigcup_{j=1}^{\ell}X_{t_j}
\]
has the property that each term $T$ of this chain is disjoint (modulo a nullset) from $\ReachOut_\Gamma(T^-)$ of the previous term $T^-$. Note that in such a setting, if $g\in L^2(\Omega)$ is such that $g$ is constant-0 almost everywhere on $T^-$, then for almost every $x\in T$, we have $(\Gamma g)(x)=\big((\Gamma\llbracket T\setminus T^-\rrbracket) g\big)(x)$. Applying this to our eigenfunction $f$ corresponding to eigenvalue $\nu$, we see if $f$ were constant-0 on $T^-$ then $\mathbbm{1}_{T\setminus T^-}f$ is either an eigenfunction of $\Gamma\llbracket T\setminus T^-\rrbracket$ corresponding to eigenvalue $\nu$ or it is constant-0. This reasoning alone yields Part~\ref{enCl:2}. For Parts~\ref{enCl:1} and~\ref{enCl:3} (which correspond to cases when $T\setminus T^-=Z_i\setminus Z_{i-1}$) we need to combine it with Proposition~\ref{prop:X0main}\ref{en:X0SpectralRadius} which asserts that $\Gamma\llbracket T\setminus T^-\rrbracket$ admits no nontrivial eigenfunction.
\end{proof}
Let us now use Claim~\ref{claimB} to conclude the statement. Suppose that $f$ is a (nonzero) eigenfunction of~$\Gamma$ corresponding to eigenvalue $\rho(\Gamma)$. Let $i\in[\ell]$ be the smallest such that $f$ is not constant-0 on $X_{t_i}$. Claim~\ref{claimB} asserts that such an $i$ exists (that is, the nonzero part of $f$ cannot be confined to $X_0$), and that $\mathbbm{1}_{X_i}f$ is an eigenfunction of $\Gamma\llbracket X_{t_i}\rrbracket$ corresponding to eigenvalue $\rho(\Gamma)$. Thus, $\rho(\Gamma\llbracket X_{t_i}\rrbracket)\ge \rho(\Gamma)$, as was needed.

\underline{The case when $I$ is infinite.} Let $I_1\subseteq I_2\subseteq \ldots \subseteq I$ be an increasing chain of finite sets that exhaust $I$. For $k\in \N$, define $\Gamma_k:=\Gamma\left\llbracket X_0\cup\bigcup_{i\in I_k} X_i\right\rrbracket$.
\begin{claim}\label{claimOp}
Digraphons $\Gamma_k$ converge to $\Gamma$ in the operator norm. 
\end{claim}
\begin{proof}[Proof of Claim~\ref{claimOp}]
Recall that the operator norm is bounded from above by the the Hilbert--Schmidt norm, and that for integral kernel operators on $L^2(\Omega)$, the latter is equal to the $L^2$-norm of the corresponding kernel. Hence, the claim follows from the obvious fact that $\|\Gamma-\Gamma_k\|_2\to 0$.
\end{proof}
We recall a well known fact that the spectrum (and hence also the spectral radius) is continuous with respect to the operator norm, when all the operators involved are compact (which digraphons are). Since by the above finite case we have $\rho(\Gamma_k)= \sup_{i\in I_k}\rho(\Gamma\llbracket X_i\rrbracket)\le \sup_{i\in I}\rho(\Gamma\llbracket X_i\rrbracket)$ for each $k\in \N$, the above continuity yields $\rho(\Gamma)\le \sup_{i\in I}\rho(\Gamma\llbracket X_i\rrbracket)$. 

\underline{The supremum is achieved.} This follows from Fact~\ref{fact:basicgraphon}\ref{en:BG3}.
\end{proof}
\begin{proof}[Proof of Proposition~\ref{prop:spectralradiusAndStrongComponents}\ref{en:CharSpectralRadius0}]
The case $I=\emptyset$ is handled by Proposition~\ref{prop:X0main}\ref{en:X0SpectralRadius}. In view of Proposition~\ref{prop:spectralradiusAndStrongComponents}\ref{en:spectralradiusMaximum}, it only needs to be proven that one digraphon $\Gamma\llbracket X_i\rrbracket$ ($i\in I$) has a nonzero eigenvalue. Of course, the choice of $i\in I$ will not be important and we will prove the assertion for every choice of $i\in I$. 

By Proposition~\ref{prop:CondensationDigraph}\ref{en:CondensationReachability} (applied on the 1-vertex walk $i$), for almost every pair $(x,y)\in X_{i}\times X_i$ there exists a positive density of paths within $X_i$ of certain length, say $\ell_{x,y}$ from $x$ to $y$. By sigma-additivity of measure, there exists a pair $(L^+,L^-)\in\N^2$ so that $\mu^2(\{(x,y)\in X_i^2:\ell_{x,y}=L^+,\ell_{y,x}=L^-\})>0$. This means that $t(C_{L^++L^-},\Gamma\llbracket X_i\rrbracket)>0$.
Applying Proposition~\ref{prop:PU} for $k=L^++L^-$, we see that $\Gamma\llbracket X_i\rrbracket$ has at least one non-zero eigenvalue.

Let us remark that a more direct proof would be using a theorem of de Pagter (see Theorem~4.2.2 in~\cite{MR1128093}), which says that non-negative, irreducible (c.f. Lemma~\ref{lem:irreducible}) compact operator has positive spectral radius.
\end{proof}

\section{Peripheral and graphical periodicity of digraphons}
In this section, we study the connection between spectral properties of a digraphon and periodicity. While the latter notion, which we formally introduce in Definition~\ref{def:graphicallyperiodic}, is new for digraphons, it is a natural counterpart to the well-known notions appearing in the theory of digraphs and even more often in the theory of Markov chains. That is, in a strongly connected digraph $H$, the following features are known to be equivalent for any given $d\in \N$ (see, e.g., Theorem~8.8.1 in~\cite{MR1829620} and Theorem~1.7 in~\cite{MR2209438}):
\begin{itemize}
    \item All directed cycles in $H$ have lengths divisible by $d$.
    \item The vertex set $V(H)$ can be partitioned into sets $P_0,P_1,\ldots,P_{d-1}$ such that each edge of $H$ goes from $P_{k-1}$ to $P_{k}$ (for some $k\in[d]$, using the cyclic notation $P_d=P_0$). The sets $P_0,\ldots,P_{d-1}$ are often called `cyclic sets'.
    \item If $\gamma$ is a largest eigenvalue in absolute value of the adjacency matrix of $H$, then the numbers $\left\{\exp(-2\pi k\ImaginaryUnit/d)\gamma\right\}_{k=1}^{d-1}$ are also eigenvalues.
    \item The entire spectrum of $H$ is symmetric (including algebraic multiplicities) under a rotation by an angle $2\pi/d$.
\end{itemize}
The counterpart of this concept for digraphons is as follows.
\begin{defi}\label{def:graphicallyperiodic}
    Suppose that $\Gamma$ is a digraphon on $\Omega$. For $d\in\N$, we say that $\Gamma$ is \emph{graphically $d$-periodic} if there exists a partition $\Omega=P_0\sqcup P_1\sqcup \ldots\sqcup P_{d-1}$ such that (using the cyclic notation $P_d=P_0$) for every $j=0,\ldots,d-1$ we have $\Gamma_{\restriction P_j\times (\Omega\setminus P_{j+1})}=0$ (almost everywhere). The sets $P_0,\ldots,P_{d-1}$ are called \emph{cyclic sets}. We say that $\Gamma$ is \emph{aperiodic} if the only integer for which this condition holds is $d=1$.
\end{defi}
A digraphon can be graphically $d$-periodic for several values of $d$. Indeed, if a digraphon is graphically $d$-periodic and $p$ divides $d$, then the partition $\Omega=Y_0\sqcup Y_1\sqcup \ldots\sqcup Y_{p-1}$, $Y_j:=\bigcup_{t:t\equiv j \mod p}P_t$ shows that it is also graphically $p$-periodic. It follows from Theorem~\ref{thm:periodicity} below that, if $\Gamma$ is strongly connected, there is a largest graphical periodicity $d$, and all the graphical periodicities are just divisors of $d$. More importantly, this number $d$ is the peripheral multiplicity of $\Gamma$.

Suppose that $S\subseteq \N$, and that $D\in\N$ and $s\in \Z$ are given. We say that $S$ is \emph{eventually $D$-periodic with shift $s$} if $S-s\subseteq D\Z$ and for some $k\in \N$ we have 
$(S-s)\cap\{k,k+1,k+2,\ldots\}=D\N\cap\{k,k+1,k+2,\ldots\}$. For example, the set $\{32,62,72,82,92,102,112,\ldots\}$ is eventually 10-periodic with shift 2, while the set $\{23,32,62,72,82,92,102,112,\ldots\}$ is not.

With these preparations, we are ready to present the main result of this section that gives a link between graphical periodicity and spectral properties. This result is accompanied by Propositions~\ref{prop:periodicblocksunique} and~\ref{prop:poweringPeriodic} that establish additional connections.
\begin{thm}\label{thm:periodicity}
Suppose that $\Gamma$ is a strongly connected digraphon on $\Omega$. Suppose that the peripheral multiplicity of $\Gamma$ is $D$. 
\begin{enumerate}[label=(\roman*)]
\item\label{en:PeriodGraphicalImpliesPeripheral} Suppose that $\Gamma$ is graphically $d$-periodic. Then the number $\exp(-2\pi \ImaginaryUnit/d) \rho(\Gamma)$ is an eigenvalue of $\Gamma$. In particular, by Theorem~\ref{thm:Schaefer74}, the number $d$ divides $D$.
\item\label{en:PeriodPeripheralImpliesGraphical} The digraphon $\Gamma$ is graphically $D$-periodic. Furthermore, there exists a partition $\Omega=P_0\sqcup P_1\sqcup \ldots\sqcup P_{D-1}$ as in Definition~\ref{def:graphicallyperiodic} such that for every $i,j\in\{0,1,\ldots,D-1\}$ and almost every pair $(x,y)\in P_i\times P_j$, the reachability sequence from $x$ to $y$ is eventually $D$-periodic with shift $j-i$.
\end{enumerate}
\end{thm}

Actually, we shall use the following corollary.
\begin{cor}\label{cor:cycles}
Suppose that $\Gamma$ is a strongly connected digraphon on $\Omega$. Suppose that the peripheral multiplicity of $\Gamma$ is $D$. Let $k\ge 3$ be arbitrary. Then $t(C_k,\Gamma)>0$ if and only if $k$ is divisible by~$D$.
\end{cor}
\begin{proof}
Theorem~\ref{thm:ActualSchaefer}\ref{en:rotationsymmetry} tells us that $\PSpec(\Gamma)$ is symmetric by rotation by an angle $2\pi/D$, including algebraic multiplicities. We can then write $\sum_{\lambda\in\PSpec(\Gamma)}\mathfrak{m}_{\Gamma}(\lambda)\cdot \lambda^k$ in $D$ different ways by applying this symmetry,
\begin{align*}
t(C_k,\Gamma)&
\eqBy{P\ref{prop:PU}}\sum_{\lambda\in\PSpec(\Gamma)}\mathfrak{m}_{\Gamma}(\lambda)\cdot\lambda^k
=\frac1D\sum_{\ell=0}^{D-1}\sum_{\lambda\in\PSpec(\Gamma)}\mathfrak{m}_{\Gamma}(\lambda)\left(\exp(-2\pi\ell\ImaginaryUnit/D)\lambda\right)^k\\
&=\frac1D\cdot \left(\sum_{\ell=0}^{D-1}\exp(-2\pi \ell k\ImaginaryUnit/D)\right)\cdot \left(\sum_{\lambda\in\PSpec(\Gamma)}\mathfrak{m}_{\Gamma}(\lambda)\cdot\lambda^k\right)
\;.    
\end{align*}
If $k$ is not divisible by $D$, $\sum_{\ell=0}^{D-1}\exp(-2\pi \ell k\ImaginaryUnit/D)=0$. On the other hand, if $k$ is divisible by $D$, then $\sum_{\ell=0}^{D-1}\exp(-2\pi \ell k\ImaginaryUnit/D)=D$. The statement follows.
\end{proof}

The proof of Theorem~\ref{thm:periodicity} has three parts. In Section~\ref{ssec:ProofEigenvalueShifted}, we prove Item~\ref{en:PeriodGraphicalImpliesPeripheral}. In Section~\ref{ssec:ProofPeriGraphi} we prove the main part of Item~\ref{en:PeriodPeripheralImpliesGraphical}, that is, the existence of a decomposition, as in Definition~\ref{def:graphicallyperiodic}. Lastly, in Section~\ref{ssec:reachability}, we prove the statement about reachability sequences.

In Proposition~\ref{prop:poweringPeriodic}, we state an easy consequence of Theorem~\ref{thm:periodicity}\ref{en:PeriodPeripheralImpliesGraphical}, namely, that the partition $\Omega=P_0\sqcup P_1\sqcup \ldots\sqcup P_{D-1}$ as in Definition~\ref{def:graphicallyperiodic} is essentially unique in the case of a strongly connected digraphon of peripheral multiplicity $D$. The proof of Proposition~\ref{prop:periodicblocksunique} is given in Section~\ref{ssec:UniquenessPeriodicity}.

\begin{prop}\label{prop:periodicblocksunique}
If $\Gamma$ is a strongly connected digraphon of peripheral multiplicity $D$, then the partition $\Omega=P_0\sqcup P_1\sqcup \ldots\sqcup P_{D-1}$ as in Definition~\ref{def:graphicallyperiodic} is unique up to a cyclic shift. That is, if $\Omega=Q_0\sqcup Q_1\sqcup \ldots\sqcup Q_{D-1}$ is as in Definition~\ref{def:graphicallyperiodic}, then for some $s\in\Z$ we have that for all $i,j\in \{0,\ldots,D-1\}$ that $P_i=_0 Q_j$ if and only if $i\equiv j+s \mod D$.
\end{prop}

The last major feature we derive in this section concerns the $D$-th power of a peripherally $D$-periodic strongly connected digraphon. This feature will be used to reduce the proof of Theorem~\ref{thm:asymptotics} only to the aperiodic case. A proof of this proposition is given in Section~\ref{ssec:prop:poweringPeriodic}.
\begin{prop}\label{prop:poweringPeriodic}
Suppose that $\Gamma$ is a strongly connected digraphon on $\Omega$. Suppose that the peripheral multiplicity of $\Gamma$ is $D$. Let $\Omega=P_0\sqcup P_1\sqcup \ldots\sqcup P_{D-1}$ be the partition provided by Theorem~\ref{thm:periodicity}\ref{en:PeriodPeripheralImpliesGraphical} for the digraphon $\Gamma$. Suppose further that $v_L$ and $v_R$ are the real eigenfunctions for the eigenvalue $\rho(\Gamma)$.

Then the decomposition of the digraphon $\Gamma^D$ into strong components (as in Theorem~\ref{thm:decompositionIntoComponents}) consists of an empty fragmented set and of strong components $P_0, P_1, \ldots, P_{D-1}$. Further, for each $i\in\{0,1,\ldots,D-1\}$, the digraphon $(\Gamma^D)_{\restriction P_i\times P_i}$ is aperiodic and we have $\rho((\Gamma^D)_{\restriction P_i\times P_i})=\frac{\rho(\Gamma)^D}{\mu(P_i)}$, and $(v_L)_{\restriction P_i}$ and $(v_R)_{\restriction P_i}$ are the left and right eigenfunctions for the eigenvalue $\rho((\Gamma^D)_{\restriction P_i\times P_i})$.
\end{prop}

\begin{rem}\label{rem:periodicitygraphons}
The issue of periodicity is particularly simple for those digraphons that are symmetric, that is, for graphons. Indeed suppose that $\Gamma$ is a connected graphon. Then Fact~\ref{fact:basicgraphon}\ref{en:BG5}. That is, either the set of peripheral eigenvalues is $\{-\rho(\Gamma),\rho(\Gamma)\}$ or $\{\rho(\Gamma)\}$. In the former case, we say that $\Gamma$ is `bipartite'. Several equivalent definitions of bipartiteness known from finite graph theory carry over to the graphon setting, see~\cite{DolezalHladky:MatchingPolytons}.
\end{rem}

\subsection{Proof of Theorem~\ref{thm:periodicity}, part~\ref{en:PeriodGraphicalImpliesPeripheral}}\label{ssec:ProofEigenvalueShifted}

Let $\Omega=X_0\sqcup X_1\sqcup \ldots\sqcup X_{d-1}$ be as in Definition~\ref{def:graphicallyperiodic}. We also use the cyclic notation $X_d=X_0$ and $X_{-1}=X_{d-1}$. Let $f\in L^2(\Omega)$ be the eigenfunction for the eigenvalue $\rho(\Gamma)$. Define a function $g\in L^2(\Omega)$ by
\[
g:=\sum_{j=0}^{d-1} \mathbbm{1}_{X_j}\cdot f \cdot \exp\left(-\frac{2\pi \ImaginaryUnit\cdot j}{d}\right) \;.
\]
In words, $g$ is obtained by taking $f$ and shifting the phase by $-\frac{2\pi j}{d}$ on every part $X_j$.

We claim $g$ is an eigenfunction for eigenvalue $\exp(-2\pi \ImaginaryUnit/d) \rho(\Gamma)$, which will prove the theorem. Indeed, let $x\in\Omega$ be arbitrary. Let $j\in\{0,1,\ldots,d-1\}$ be such that $x\in X_j$. Recall that by the key property of Definition~\ref{def:graphicallyperiodic}, except for a nullset of exceptional $x\in X_j$, we have that 
\begin{equation}\label{eq:MONull}
\mbox{$\Gamma(\cdot,x)$ is zero almost-everywhere on $\Omega\setminus X_{j-1}$}\;.
\end{equation}
Then
\begin{align*}
(\Gamma g)(x)&=
\int_{z\in \Omega}\Gamma(z,x)g(z)\diff z
\eqByRef{eq:MONull}
\int_{z\in X_{j-1}}\Gamma(z,x)g(z)\diff z\\
&
=
\int_{z\in X_{j-1}}\Gamma(z,x)\cdot f(z)\cdot\exp\left(-\frac{2\pi \ImaginaryUnit\cdot (j-1)}{d}\right)\diff z\\
&=
\exp\left(-\frac{2\pi \ImaginaryUnit\cdot (j-1)}{d}\right)\cdot \int_{z\in X_{j-1}}\Gamma(z,x)\cdot f(z)\diff z\\
&
\eqByRef{eq:MONull}
\exp\left(-\frac{2\pi \ImaginaryUnit\cdot (j-1)}{d}\right)\cdot \int_{z\in \Omega}\Gamma(z,x)\cdot f(z)\diff z\\
\JUSTIFY{$f$ eigenfunction}
&=\exp\left(-\frac{2\pi \ImaginaryUnit\cdot (j-1)}{d}\right) \rho(\Gamma) \cdot f(x)=
\exp\left(-\frac{2\pi \ImaginaryUnit\cdot (j-1)}{d}\right) \rho(\Gamma) \cdot g(x) \exp\left(\frac{2\pi \ImaginaryUnit\cdot j}{d}\right)\\
&=\exp\left(\frac{2\pi \ImaginaryUnit}{d}\right)\rho(\Gamma)  g(x)\;,
\end{align*}
as was needed.

\subsection{Proof of Theorem~\ref{thm:periodicity}, the main part of~\ref{en:PeriodPeripheralImpliesGraphical}}\label{ssec:ProofPeriGraphi}
We will use a trivial inequality for complex numbers. Recall that the angular component of the polar coordinates of a complex number is called the \emph{phase}. The phase is uniquely defined except for number~0 which can have every phase. 
\begin{fact}
\label{fact:ComplexNumbers}
Suppose that $\alpha$ is a finite measure on a measure space $\Xi$, $f:\Xi\rightarrow\mathbb{C}$ is a bounded function, and we have $\left|\int f\mathsf{d}\alpha\right|=\int\left|f\right|\mathsf{d}\alpha$. Then the phase of $f(x)$ is constant $\alpha$-almost everywhere.
\end{fact}

Let $f$ be an eigenfunction corresponding to eigenvalue $\tau:=\exp(-\frac{2\pi \ImaginaryUnit}{D})\rho(\Gamma)$. That is, for almost every $x\in \Omega$, we have
\begin{equation}\label{eq:nok}
\int_{y}\Gamma(y,x) f(y)\diff\mu(y)=\tau f(x)\;. 
\end{equation}
Define a function $g:\Omega\to\R$, $g(x):=\left|f(x)\right|$. Take an arbitrary $x\in \Omega$. We have
\begin{equation}\label{eq:littletired}
\left(\Gamma g\right)(x)=\int_{y}\Gamma(y,x)g(y)\diff\mu(y)\ge\left|\int_{y}\Gamma(y,x)f(y)\diff\mu(y)\right|=\left|\tau f(x)\right|=\rho(\Gamma) g(x)\;.    
\end{equation}
So, Theorem~\ref{thm:Schaefer74}\ref{en:PrincipalEigevectorInequality} tells us that $g$ is an
eigenfunction for eigenvalue $\rho(\Gamma)$, that is, $\Gamma g=\rho(\Gamma)g$. That is, for almost every $x\in X$, there is an equality between the second and the third term of~\eqref{eq:littletired},
\begin{equation}\label{eq:kon}
\int_{y}\Gamma(y,x)\left|f(y)\right|\diff\mu(y)=\left|\int_{y}\Gamma(y,x)f(y)\diff\mu(y)\right|\;. 
\end{equation}
The next claim it crucial. It asserts that for a given $x$, $f(x)$ differs in phase by $-\frac{2\pi}{D}$ from $f(y)$ for almost all inneighbors $y$ of $x$.
\begin{claim}\label{claimC}
    Suppose that $x\in\Omega$ is given and that it satisfies~\eqref{eq:nok} and~\eqref{eq:kon}. Then for almost all $y\in\Omega$ with $\Gamma(y,x)>0$, the number $f(y)$ has the phase of $f(x)$ shifted by $-\frac{2\pi}{D}$.
\end{claim}
\begin{proof}[Proof of Claim~\ref{claimC}]
Define a measure $\alpha$ on $\Omega$ by $\alpha(T):=\int_{y\in T}\Gamma(y,x)\diff y$. We can rewrite~\eqref{eq:kon} as $\left|\int f\mathsf{d}\alpha\right|=\int\left|f\right|\mathsf{d}\alpha$. In particular, Fact~\ref{fact:ComplexNumbers} tells us that for $\alpha$-almost all $y$'s, $f(y)$ has the same phase (note that this does not say anything about $y$'s with $\Gamma(y,x)=0$). 

The rest follows from~\eqref{eq:nok}. Indeed, the right-hand side has the phase of $f(x)$ shifted by $-\frac{2\pi}{D}$. As for the left-hand side, we already proved that all contributing terms $\Gamma(y,x)f(y)$ have the same phase (recall that $\Gamma(y,x)>0$, and this factor does not change the phase), say $\theta$. It follows that this phase $\theta$ has to be the phase appearing on the right-hand side.
\end{proof}
Fix $x$ that satisfies~\eqref{eq:nok} and~\eqref{eq:kon}. It follows by induction on $\ell=1,2,\ldots$ that for almost all $y$ with $x$ reachable from $y$ in $\ell$ steps, the values $f(y)$ has phase of $f(x)$ shifted $-\frac{2\pi\ell}{D}$ (these shifts are $D$-periodic). Theorem~\ref{thm:Schaefer74}\ref{en:LeadingEigenvectorPositive} tells us that $g$ is strictly positive, almost everywhere. In other words, $f$ is not zero almost everywhere, and hence the phase of $f(y)$ is unique for almost all $y$. Since, for almost every $y$, $x$ can be reached from $y$ in a finite number of steps (recall Proposition~\ref{prop:reachabilityFromAnywhereToAnywhere}), we conclude that for almost all $y\in\Omega$, the phase of $f(y)$ differs from the phase of $f(x)$ by $-\frac{2\pi\ell}{D}$ for some $\ell=0,1,\ldots,D-1$. We partition almost all elements of $\Omega$ according to this phase shift, $\Omega =_0 P_0\sqcup\ldots\sqcup P_{D-1}$. We now need to verify that this partition satisfies Definition~\ref{def:graphicallyperiodic}. Pick $j\in\{0,1\ldots,D-1\}$. It is our task to show that for almost all $\tilde{x}\in P_{j+1}$ we have that almost all the inneighbors of $\tilde{x}$ lie in $P_j$. But this is precisely what Claim~\ref{claimC} tells us.

\subsection{Proof of Theorem~\ref{thm:periodicity}, reachability sequences}\label{ssec:reachability}
We use the partition $\Omega =_0 P_0\sqcup\ldots\sqcup P_{D-1}$ obtained in Section~\ref{ssec:ProofPeriGraphi}
The proof has two parts. In the first (and easy) part, we prove that for almost every $(x,y)\in P_i\times P_j$ we have that if $k\in \N$ is such that $k-(j-i)$ is not divisible by $D$, then $\Gamma^k(x,y)=0$. In the second (and longer) part, we will prove that eventually all reachability lengths of one congruence class modulo $D$ occur between $x$ and $y$.

\subsubsection{First part: reachability lengths of wrong congruences do not appear}\label{sssec:POI}
Equivalently, the first part amounts to proving that $\int_{(x,y)\in P_i\times P_j}\Gamma(x,y)=0$.

We have
\[
\int_{(x,y)\in P_i\times P_j}=\int_{(x_0,x_1,\ldots,x_k)\in P_i\times \Omega^{k-1}\times P_j}\prod_{\ell=1}^k\Gamma(x_{\ell-1},x_\ell)
\;.
\]
We split the domain $P_i\times \Omega^{k-1}\times P_j$ into $2^{k-1}$ parts, $P_i\times \Omega^{k-1}\times P_j=P_i\times \left(\bigcup_{a_1,\ldots,a_{k-1}\in\{\oplus,\ominus\}}\prod_{\ell=1}^{k-1} A_\ell^{a_\ell}\right)\times P_j$, where $A^\oplus_\ell:=P_{i+\ell}$ (with cyclic notation of indices) and $A^\ominus_\ell=\Omega\setminus P_{i+\ell}$. Definition~\ref{def:graphicallyperiodic} tells us that $\Gamma$ is constant-0 on $P_i\times A_1^\ominus$ and also $A_{\ell-1}^\oplus\times A_\ell^\ominus$ on for every $\ell\in[k]$. We also have that $\Gamma$ is constant-0 on $A_k^\oplus \times P_j$. Put together, we conclude that for every choice of $a_1,\ldots,a_{k-1}\in\{\oplus,\ominus\}$, we have $\int_{(x_0,x_1,\ldots,x_k)\in P_i\times \prod_{\ell=1}^{k-1} A_\ell^{a_\ell}\times P_j}\prod_{\ell=1}^k\Gamma(x_{\ell-1},x_\ell)=0$, as was needed for the first part.

\subsubsection{Second part: reachability lengths of the right congruences eventually all appear}
By Corollary~\ref{cor:cycles}, we have $t(C_D,\Gamma)>0$. For each $i\in\{0,1,\ldots,D-1\}$, let $P_i^*\subseteq P_i$ consist of those $x$ with $t^{\bullet}_{x}(C_k^{\bullet},\Gamma)>0$. By Fubini's theorem, the sets $P_0^*,P_1^*,\ldots,P_{D-1}^*$ have positive measures.

Now, suppose that $i$ and $j$ are given and that $x\in P_i$ and $y\in P_j$. Find a number $\ell\in\N$ so that $\int_{z\in P_0^*}t_{x,z}^{\bullet\bullet}(P_{\ell},\Gamma)>0$. Such a number exists by Proposition~\ref{prop:reachabilityFromAnywhereToAnywhere}. Let $P_0^{**}\subseteq P_0^*$ be the set of $z$'s such that $t_{x,z}^{\bullet\bullet}(P_{\ell},\Gamma)>0$. Hence, $P_0^{**}$ is of positive measure. Find a number $h\in\N$ so that $\int_{z\in P_0^{**}}t_{z,y}^{\bullet\bullet}(P_{h},\Gamma)>0$. Such a number exists by Proposition~\ref{prop:reachabilityFromAnywhereToAnywhere}.
\begin{claim}\label{cl:reachabilityfinished}
    For every $j\in \N_0$, we have $t_{x,y}^{\bullet\bullet}(P_{\ell+h+jD},\Gamma)>0$.
\end{claim}
Before proving Claim~\ref{cl:reachabilityfinished}, let us show that this implies the statement. Indeed, Claim~\ref{cl:reachabilityfinished} asserts that the reachability sequence from $x$ to $y$ is eventually $D$-periodic with some shift. But Section~\ref{sssec:POI} tells us that the only available shift is $j-i$.
\begin{proof}[Proof of Claim~\ref{cl:reachabilityfinished}]
Indeed, we have a positive density of paths from $x$ to $P^{**}_0$. From each element of $P^{**}_0$, wrap $j$ cycles of length $D$ around that vertex, and finally, use paths of length from $P^{**}_0$ to $y$. 
\end{proof}

\subsection{Proof of Proposition~\ref{prop:periodicblocksunique}}\label{ssec:UniquenessPeriodicity}
We need to prove that if $\Omega=P_0\sqcup P_1\sqcup \ldots\sqcup P_{D-1}$ provided by Theorem~\ref{thm:periodicity}\ref{en:PeriodPeripheralImpliesGraphical}, and $\Omega=Q_0\sqcup Q_1\sqcup \ldots\sqcup Q_{D-1}$ is another partition of $\Omega$ which is not a cyclic shift of $P_0\sqcup P_1\sqcup \ldots\sqcup P_{D-1}$. We need to prove that $\Omega=Q_0\sqcup Q_1\sqcup \ldots\sqcup Q_{D-1}$ does not satisfy the properties of Definition~\ref{def:graphicallyperiodic}. Below, we formalize the assumption that $Q_0\sqcup Q_1\sqcup \ldots\sqcup Q_{D-1}$ is not a cyclic shift of $P_0\sqcup P_1\sqcup \ldots\sqcup P_{D-1}$ but here we mention that the reason could be twofold: (i)~there exist $i,j,k$ with $j\neq k$ so that both $P_i\cap Q_j$ and $P_i\cap Q_k$ have positive measures, or~(ii)~the unordered partitions $\{P_0,\ldots,P_{D-1}\}$ and $\{Q_0,\ldots,Q_{D-1}\}$ are the same (modulo nullsets) but the bijection between these two partitions is not a cyclic shift.

Since $Q_0\sqcup Q_1\sqcup \ldots\sqcup Q_{D-1}$ is not a cyclic shift of $P_0\sqcup P_1\sqcup \ldots\sqcup P_{D-1}$ then there exist indices (not necessarily distinct) $i,j,h,k\in\{0,\ldots,D-1\}$ with 
\begin{equation}\label{someprvni}
j-i\not\equiv h-k\mod D    
\end{equation}
and the property that both $S:=P_i\cap Q_k$ and $T:=P_j\cap Q_h$ have positive measures. Then on the one hand, Theorem~\ref{thm:periodicity}\ref{en:PeriodPeripheralImpliesGraphical} tells us that for almost every $(x,y)\in S\times T$, the reachability sequence from $x$ to $y$ is eventually $D$-periodic with shift $j-i$. In particular, for some 
\begin{equation}\label{someL}
L\equiv j-i\mod D    
\end{equation}
we have 
\begin{equation}\label{osm}
\int_{(x,y)\in S\times T}t^{\bullet\bullet}_{x,y}\left(P_{L}^{\bullet\bullet},\Gamma\right)>0\;.
\end{equation}

On the other hand, the claim below shows that the assumption that $Q_0\sqcup Q_1\sqcup \ldots\sqcup Q_{D-1}$ satisfies Definition~\ref{def:graphicallyperiodic} would lead to a contradiction.
\begin{claim}\label{cl:strankonaskoc}
If $Q_0\sqcup Q_1\sqcup \ldots\sqcup Q_{D-1}$ satisfied Definition~\ref{def:graphicallyperiodic}, and for a number $r\in \N$ and a sequence $k=:t_1,t_2,\ldots,t_r\in\{0,\ldots,D-1\}$ we had that $\int_{x_1\in Q_{t_1}}\int_{x_2\in Q_{t_2}}\cdots \int_{x_r\in Q_{t_r}}\prod_{q=1}^{r-1}\Gamma(x_q,x_{q+1})>0$, then we would have that $t_2\equiv k+1\mod D$, $t_{r}\equiv k+r-1\mod D$.
\end{claim}
Before we prove Claim~\ref{cl:strankonaskoc}, let us show that the assumption that $Q_0\sqcup Q_1\sqcup \ldots\sqcup Q_{D-1}$ satisfied Definition~\ref{def:graphicallyperiodic} would lead to a contradiction. Indeed, in that case we would have
\begin{align*}
\int_{(x,y)\in S\times T}t^{\bullet\bullet}_{x,y}\left(P_{L}^{\bullet\bullet},\Gamma\right)&\le 
\int_{x_1\in Q_k}\int_{x_2\in \Omega}\cdots \int_{x_{L}\in \Omega}\int_{x_{L+1}\in Q_h}\prod_{q=1}^{L}\Gamma(x_q,x_{q+1})
\\
&
=
\sum_{t_2=0}^{D-1}\sum_{t_3=0}^{D-1}\cdots \sum_{t_L=0}^{D-1}
\int_{x_1\in Q_k}\int_{x_2\in Q_{t_2}}\cdots \int_{x_{L}\in Q_{t_L}}\int_{x_{L+1}\in Q_h}\prod_{q=1}^{L}\Gamma(x_q,x_{q+1})
\;.    
\end{align*}
Note that~\eqref{someprvni} and~\eqref{someL} imply that $h\not\equiv k+L\mod D$. Thus, Claim~\ref{cl:strankonaskoc} yields that for every choice of $t_2,\ldots,t_L$ we have $\int_{x_1\in Q_k}\int_{x_2\in Q_{t_2}}\cdots \int_{x_{L}\in Q_{t_L}}\int_{x_{L+1}\in Q_h}\prod_{q=1}^{L}\Gamma(x_q,x_{q+1})=0$. Hence, $\int_{(x,y)\in S\times T}t^{\bullet\bullet}_{x,y}\left(P_{L}^{\bullet\bullet},\Gamma\right)=0$, a contradiction to~\eqref{osm}.
\begin{proof}[Proof of Claim~\ref{cl:strankonaskoc}]
This easily follows by induction on $r=1,2,\ldots$ and the key property of Definition~\ref{def:graphicallyperiodic}.
\end{proof}

\subsection{Proof of Proposition~\ref{prop:poweringPeriodic}}\label{ssec:prop:poweringPeriodic}
The proof consists of two complementing parts concerning strong connectivity and aperiodicity (namely an argument about disconnectedness in Lemma~\ref{lem:powerDisconnected} and an argument about connectedness in Lemma~\ref{lem:powerWellConnected}), and Lemma~\ref{lem:Prop65Spectral} which covers the spectral part of the statement.
\begin{lem}\label{lem:powerDisconnected}
In the setting of Proposition~\ref{prop:poweringPeriodic}, we have for each $i\in\{0,1,\ldots,D-1\}$ that $\int_{P_i\times (\Omega\setminus P_i)}\Gamma^D=\int_{(\Omega\setminus P_i)\times P_i}\Gamma^D=0$.
\end{lem}
\begin{lem}\label{lem:powerWellConnected}
In the setting of Proposition~\ref{prop:poweringPeriodic}, we have for each $i\in\{0,1,\ldots,D-1\}$ and for almost every pair $(x,y)\in P_i\times P_i$ that the reachability sequence from $x$ to $y$ in the digraphon $\Gamma^D_{\restriction P_i\times P_i}$ is eventually 1-periodic.
\end{lem}
\begin{lem}\label{lem:Prop65Spectral}
In the setting of Proposition~\ref{prop:poweringPeriodic}, we have for each $i\in\{0,1,\ldots,D-1\}$ that $\rho((\Gamma^D)_{\restriction P_i\times P_i})=\frac{\rho(\Gamma)^D}{\mu(P_i)}$, and $(v_L)_{\restriction P_i}$ and $(v_R)_{\restriction P_i}$ are the left and right eigenfunctions for the eigenvalue $\rho((\Gamma^D)_{\restriction P_i\times P_i})$.
\end{lem}
It is easy to see that Proposition~\ref{prop:poweringPeriodic} follows. Indeed the fact implied by Lemma~\ref{lem:powerWellConnected} that between almost every pair $(x,y)\in P_i\times P_i$, $y$ is reachable from $x$ in $\Gamma^D$ implies that $P_i$ is strongly connected. Lemma~\ref{lem:powerDisconnected}, on the other hand, shows that $P_i$ also satisfies Definition~\ref{def:component}\ref{en:StrongComp}.
It remains to argue that $\Gamma^D_{\restriction P_i\times P_i}$ is aperiodic. Indeed, if $\Gamma^D_{\restriction P_i\times P_i}$ were $d$-periodic for some $d>1$, then by Theorem~\ref{thm:periodicity}, for almost every $(x,y)\in P_i\times P_i$, the reachability sequence from $x$ to $y$ would be eventually $d$-periodic (with a certain shift $s_{x,y}$), a contradiction with the eventual 1-peridicity given by Lemma~\ref{lem:powerWellConnected}. The spectral part of Proposition~\ref{prop:poweringPeriodic} is covered by Lemma~\ref{lem:Prop65Spectral}. Hence, Proposition~\ref{prop:poweringPeriodic} is proven.

Crucial to our proofs of Lemmas~\ref{lem:powerDisconnected} and~\ref{lem:powerWellConnected} will be the following two simple observations. The first one follows directly from Definition~\ref{def:digraphonpower}: For every $k\in\N$ and $(x,y)\in \Omega^2$, we have
\begin{equation}\label{eq:pathsinpowers}
t^{\bullet\bullet}_{x,y}(P_{k}^{\bullet\bullet},\Gamma^D)=t^{\bullet\bullet}_{x,y}(P_{Dk}^{\bullet\bullet},\Gamma)\;.
\end{equation}
The second one combines~\eqref{eq:pathsinpowers} with Theorem~\ref{thm:periodicity}. Namely, the furthermore part of Theorem~\ref{thm:periodicity}\ref{en:PeriodPeripheralImpliesGraphical} tells us that if $x\in P_i$ then expressing $t^{\bullet\bullet}_{x,y}(P_{k}^{\bullet\bullet},\Gamma^D)$ as~\eqref{eq:densPathTwoTerminals}, the integrand is nonzero only if $x_2,\ldots,x_k,y\in P_i$. Thus, for $(x,y)\in P_i\times P_i$,
\begin{equation}\label{eq:pathsinpowersbetter}
t^{\bullet\bullet}_{x,y}(P_{k}^{\bullet\bullet},\Gamma^D_{\restriction P_i\times P_i})=(\frac1{\mu(P_i)})^k \cdot t^{\bullet\bullet}_{x,y}(P_{Dk}^{\bullet\bullet},\Gamma)\;.
\end{equation}
where the term $(\frac1{\mu(P_i)})^k$ comes from the renormalization of the measure as in Definition~\ref{def:restriction}.

\begin{proof}[Proof of Lemma~\ref{lem:powerDisconnected}]
Equivalently, we need to prove that for almost every $(x,y)\in P_i\times(\Omega\setminus P_i)$ we have $t^{\bullet\bullet}_{x,y}(P_{1}^{\bullet\bullet},\Gamma^D_{\restriction P_i\times P_i})=t^{\bullet\bullet}_{y,x}(P_{1}^{\bullet\bullet},\Gamma^D_{\restriction P_i\times P_i})=0$. Thanks to~\eqref{eq:pathsinpowers}, this is equivalent to asking whether the reachability sequence from $x$ to $y$ and from $y$ to $x$ in the digraphon $\Gamma$ exclude the number $D$. Theorem~\ref{thm:periodicity}\ref{en:PeriodPeripheralImpliesGraphical} tells us that this is indeed the case (since the shift is nonzero).
\end{proof}
\begin{proof}[Proof of Lemma~\ref{lem:powerWellConnected}]
Theorem~\ref{thm:periodicity}\ref{en:PeriodPeripheralImpliesGraphical} tells us that for almost every $(x,y)\in P_i\times P_i$, the reachability sequence from $x$ to $y$ in the digraphon $\Gamma$ is eventually $D$-periodic with zero shift. We use~\eqref{eq:pathsinpowersbetter}, and see that the reachability sequence from $x$ to $y$ in the digraphon $\Gamma^D_{\restriction P_i\times P_i}$ is eventually $1$-periodic. This also implies the statement about strong connectedness.
\end{proof}

\begin{proof}[Proof of Lemma~\ref{lem:Prop65Spectral}]
The claim that $\rho((\Gamma^D)_{\restriction P_i\times P_i})=\frac{\rho(\Gamma)^D}{\mu(P_i)}$ is equivalent to $\rho((\Gamma^D\llbracket P_i\rrbracket)=\rho(\Gamma)^D$. The $\le$-inequality follows from Proposition~\ref{prop:spectralradiusAndStrongComponents}. For the $\ge$-inequality, it suffices to show that $(v_L)_{\restriction P_i}$ and $(v_R)_{\restriction P_i}$ are the left and right eigenfunctions for the eigenvalue $\rho((\Gamma^D)_{\restriction P_i\times P_i})$. We show the argument for $(v_L)_{\restriction P_i}$; the argument for $(v_R)_{\restriction P_i}$ is verbatim. This claim amounts to proving that for the function $f\in L^2(\Omega)$, $f:=v_L \cdot \mathbbm{1}_{P_i}$ we have $\Gamma^D f=\rho(\Gamma)^D\cdot f$. This is obvious.
\end{proof}

\section{Asymptotics for values of high powers of a digraphon}
In this section, we present one of the main theorems of the paper. It asserts that high powers $\Gamma^k$ of a strongly connected digraphon $\Gamma$ can be approximately expressed in terms of the spectral radius, and the corresponding left and right eigenfunctions. The theorem is a functional-analytic counterpart to well-known results from the Perron--Frobenius theory. Namely, if $A$ is a square matrix with positive entries and $k$ is large, then $A^k\approx \lambda^k v_{right}v_{left}$, where $\lambda$ is the largest eigenvalue of $A$, $v_{left}$ and $v_{right}$ are the unique left and right eigenvectors corresponding to it, and the product $v_{right}v_{left}$ is a matrix of the same dimensions as $A$. See for example Theorem~8.2.8(f) in~\cite{MatrixAnalysis}.
\begin{thm}\label{thm:asymptotics}
Suppose that $\Gamma$ is a strongly connected digraphon on ground set $\Omega$.
We assume that there are left and right real eigenfunctions
$v_L$, $v_R$ for the eigenvalue $\rho(\Gamma)$ satisfying $\langle v_L, v_R \rangle = 1$.

Let the peripheral multiplicity of $\Gamma$ be $D$. Suppose that $\Omega=P_0\sqcup P_1\sqcup \ldots\sqcup P_{D-1}$ is a decomposition as in Theorem~\ref{thm:periodicity}\ref{en:PeriodPeripheralImpliesGraphical}.

Let $\rho:=\rho(\Gamma)$. There exists a number $\alpha\in (0,\rho)$ with the following property: For every $i,j\in\{0,\ldots,D-1\}$ and almost every $(x,y)\in P_i\times P_j$ we have
\begin{equation}\label{TE}
\Gamma^\ell(x,y) = \begin{cases}
                    \rho^\ell v_R(x) v_L(y) + O(\alpha^\ell)&\quad\mbox{if $\ell\equiv j-i \mod D$, or}  \\
                    0&\quad\mbox{otherwise,} 
                   \end{cases} 
\end{equation}
as $\ell\to\infty$. The term in $O(\cdot )$ does not depend on $x$ and $y$.
\end{thm}

The proof of Theorem~\ref{thm:asymptotics} has two parts. First, we prove it in the aperiodic case, as stated in Proposition~\ref{prop:asymptoticAperiodic}. We give a proof of Proposition~\ref{prop:asymptoticAperiodic} in Section~\ref{ssec:proofPropAsymptoticAperiodic}. Then, in Section~\ref{ssec:asymptoticAperiodicPeriodic}, we show that the aperiodic case implies the full Theorem~\ref{thm:asymptotics}.
\begin{prop}\label{prop:asymptoticAperiodic}
Suppose that $W$ is a strongly connected aperiodic digraphon on ground set $\Omega$. Let $\tau:=\rho(W)$. We assume that there are left and right real eigenfunctions
$w_L$, $w_R$ for the eigenvalue $\rho(\Gamma)$ satisfying $\langle w_L, w_R \rangle = 1$. There exists a number $\beta\in (0,\tau)$ such that for almost every $(x,y)\in \Omega^2$ we have $W^\ell(x,y) = \tau^\ell w_R(x) w_L(y) + O(\beta^\ell)$ as $\ell\to\infty$. The term in $O(\cdot )$ does not depend on $x$ and $y$.
\end{prop}

\begin{rem}\label{rem:asymptoticsforgraphons}
Theorem~\ref{thm:asymptotics} is much easier and has been widely used when $\Gamma$ is a graphon. Indeed, in that case one can use the spectral decomposition~\eqref{eq:spectraldecomposition}, 
$\Gamma(x, y) = \sum_{i} \lambda_i f_i(x) f_i(y)$, where $\{f_i(x)\}_i$ are orthonormal eigenfunctions of $\Gamma$ (the symmetry of $\Gamma$ gives that each eigenfunction is left and right at the same time), and $\{\lambda_i\}_i$ are the corresponding eigenvalues. Orthonormality means that for each $k\in \N$, $\Gamma^k(x,y)=\sum_{i} \lambda_i^k f_i(x) f_i(y)$. That is, we have an explicit description of the error term in~\eqref{TE}. Recall that Remark~\ref{rem:periodicitygraphons} tells us that the issue of periodicity is particularly simple for graphons.
\end{rem}

\subsection{Proof of Proposition~\ref{prop:asymptoticAperiodic}}\label{ssec:proofPropAsymptoticAperiodic}

Let $T$ be the integral kernel operator corresponding to the digraphon $W$. We will decompose $T$ as $T=T_1+T_2$, where $T_1$ will be a rank-1 operator with spectrum $\Spec(T_1)=\{\tau\}$ and $T_2$ will satisfy $\rho(T_2)<\tau$. We will use tools from Section~\ref{ssec:ReducingPairs} to this end. Firstly, note that Fact~\ref{fact:basicgraphon}\ref{en:BG3} implies that $\{\tau\}$ is indeed a spectral set of $T$. Thus, we can apply Theorem~\ref{thm:reducing}. Let $L^2(\Omega)=X_1\oplus X_2$ be the corresponding reducing pair, and let $P_1$ and $P_2$ be the corresponding projections (given by Proposition~\ref{prop:reducingpair}). Define $T_1:=T\circ P_1$ and $T_2:=T\circ P_2$. We indeed have 
\begin{equation}\label{eq:Tsum}
    T=T\circ \Identity\eqBy{Fact~\ref{fact:projectionsId}} T\circ (P_1+P_2)=(T\circ P_1)+ (T\circ P_2)= T_1+T_2\;.
\end{equation}

\begin{lem}\label{lem:proofReduced}
In the setting above, we have the following properties.
\begin{enumerate}[label=(\roman*)]
    \item\label{r1} $\dim(X_1)=1$,
    \item\label{r2} $\rho(T_2)<\tau$.
\end{enumerate}
\end{lem}
\begin{proof}[Proof of Lemma~\ref{lem:proofReduced}\ref{r1}]
By Theorem~\ref{thm:Schaefer74}\ref{en:EigenvaluesAndSimple}, $\tau$ has algebraic multiplicity~1 for operator $T$. Thus, the claim follows from Lemma~\ref{lem:ReducingPairAM}.
\end{proof}
\begin{proof}[Proof of Lemma~\ref{lem:proofReduced}\ref{r2}]
It follows by Theorem~\ref{thm:reducing} that $\Spec(T_2)=\Spec(T)\setminus \{\tau\}$. Recall that $W$ is aperiodic. In particular, by Theorem~\ref{thm:periodicity}, $\tau$ is the only peripheral eigenvalue of $W$. Combined with the fact that the only possible accumulation point of $\Spec(T)$ is~0, we get that $\rho(T_2)<\tau$.
\end{proof}
Set $\beta:=(\tau+\rho(T_2))/2$. We have $\beta\in (0,\tau)$, as required in the statement of the proposition.

The next lemma shows that $T^k = T_1^k + T_2^k$ for every $k\in \N$.
\begin{lem}\label{lem:Tk}
    For every $k\in \N$, we have $T^k = T_1^k + T_2^k$.
\end{lem}
\begin{proof}
First, let us show that $T_1 \circ T_2 = 0$. Let $y\in L^2(\Omega)$. Since $L^2(\Omega)=X_1\oplus X_2$, we can write $y=y_1+y_2$, where $y_i\in X_i$. Observe that by Proposition~\ref{prop:reducingpair}, for $i=1,2$, $P_{3-i}(y_i)=0$, and also $T\circ P_i=P_i \circ T$. Thus,
\[
(T_1 \circ T_2)(y)=(T\circ P_1\circ T\circ P_2)(y_1+y_2)=(T\circ P_1\circ T)(y_2)=(T\circ T\circ P_1)(y_2)=(T\circ T)(0)=0
\]
as was needed.
The same method shows that $T_2 \circ T_1 = 0$. 

To conclude, in the expansion of $T^k = (T_1+T_2)^k$, all the terms except $T_1^k$ and $T_2^k$ vanish. The statement of the lemma follows.
\end{proof}

Recall that we are given real left and right eigenfunctions
$w_L$ and $w_R$ for the eigenvalue $\tau$ satisfying $\langle w_L, w_R \rangle = 1$.
Lemma~\ref{lem:eigenvectorsbounded} (applied either to $W$ or to $W^\top$) tells us that 
\begin{equation}\label{eq:wlr}
\|w_L\|_\infty<\infty  \quad\text{and}\quad \|w_R\|_\infty<\infty.
\end{equation}

The next lemma describes $T_1$ and $T_2$ as integral kernel operators.
\begin{lem}\label{lem:describeTiAsWi}$~$
\begin{enumerate}[label=(\roman*)]
    \item\label{en:W1} $T_1$ is an integral kernel operator with kernel $W_1$, $W_1(x,y)=\tau w_R(x) w_L(y)$. Furthermore, for each $k\in \N$, we have $W_1^k(x,y)=\tau^k w_R(x) w_L(y)$.
    \item\label{en:W2} $T_2$ is an integral kernel operator with kernel $W_2\in L^\infty(\Omega^2)$.
\end{enumerate}
\end{lem}
\begin{proof}[Proof of Proposition~\ref{lem:describeTiAsWi}\ref{en:W1}]
Clearly, $w_L \in X_1$. Since $\dim(X_1)=1$ (Lemma~\ref{lem:proofReduced}\ref{r1}), the function $w_L$ spans $X_1$. For every $f \in L^2(\Omega)$, we have
$$
\langle \tau f - T(f), w_R \rangle =
\langle \tau f , w_R \rangle -\langle T(f), w_R \rangle \eqByRef{eq:Adjoint}
\tau \langle f, w_R \rangle - \langle f, T^*(w_R) \rangle = 
\tau \langle f, w_R \rangle - \tau \langle f,w_R \rangle = 0 \;.
$$
This shows that $w_R\in \image (\tau \Identity - T)^\bot=X_2^\bot$, where the last step used Lemma~\ref{lem:ReducingPairAM}. Thus,
$$
P_1(f) = \langle f, w_R \rangle w_L\;,
$$
and consequently,
$$
T_1(f) = \rho(T) \langle f, w_R \rangle w_L \;.
$$
By induction, we obtain
$$
T_1^k(f) = \rho(T)^k \langle f, w_R \rangle w_L \;.
$$
So, we have
\begin{equation*}
T_1^k(f)(x) = \rho(T)^k w_L(x) \int_\Omega f(y) w_R(y) \diff \mu(y) \;,
\end{equation*}
as claimed.

\end{proof}
\begin{proof}[Proof of Proposition~\ref{lem:describeTiAsWi}\ref{en:W2}]
We have $T_2=T-T_1$ by~\eqref{eq:Tsum}. Hence, $T_2$ corresponds to an integral kernel operator $W_2:=W-W_1$. $W$ is a digraphon, and hence bounded in the $L^\infty$-norm. $W_1$ has finite $L^\infty$-norm, since by Proposition~\ref{lem:describeTiAsWi}\ref{en:W1} it can be expressed as a pointwise product of $w_L$ and $w_R$, which satisfy~\eqref{eq:wlr}.
\end{proof}

We are now in a position to express $W^k$. Lemma~\ref{lem:Tk} tells us that $W^k=W_1^k+W_2^k$. Lemma~\ref{lem:describeTiAsWi}\ref{en:W1} tells us that $W_1^k(x,y)=\tau^k w_R(x) w_L(y)$. We now turn to $W_2^k$. Recall that $\|W_2\|_\infty <\infty$ by Lemma~\ref{lem:describeTiAsWi}\ref{en:W2}. Proposition~\ref{prop:asymptotics} applied on $W_2$ and with $\beta$ asserts that $\|W_2^k\|_\infty=O(\beta^k)$ as $k\to\infty$. The statement of the proposition follows.

\subsection{Deducing Theorem~\ref{thm:asymptotics} from Proposition~\ref{prop:asymptoticAperiodic}}\label{ssec:asymptoticAperiodicPeriodic}

Suppose the setting of Theorem~\ref{thm:asymptotics}. Rather than proving~\eqref{TE} for a universal constant $\alpha<\rho$, we will find one number  $\alpha_{i,j}\in (0,\rho)$ for each $i,j\in\{0,\ldots,D-1\}$ so that~\eqref{TE} holds with error term $O(\alpha_{i,j}^\ell)$. This is sufficient, as then we can take $\alpha:=\max_{i,j} \alpha_{i,j}$. The `otherwise' part of~\eqref{TE} is trivial. So, $\ell=MD+j-i$, for large $M\in \N$.

\underline{First, we prove the statement when $i=j$.} In that case, for each $(x,y)\in P_i\times P_i$,
\begin{equation}\label{eq:SubstIntoMe}
\Gamma^\ell(x,y)=\mu(P_i)^{M-1}\cdot\left((\Gamma^D)_{\restriction P_i\times P_i}\right)^M(x,y)\;,     
\end{equation}
where the factor $\mu(P_i)^{M-1}$ comes from the rescaling of the probability space $\mu_{P_i}$ as in Definition~\ref{def:restriction}. Define $W:=\Gamma^D$. Let $\tau:=\rho(W)$. By Proposition~\ref{prop:poweringPeriodic}, for we have that $(\Gamma^D)_{\restriction P_i\times P_i}$ is strongly connected, aperiodic, $\rho((\Gamma^D)_{\restriction P_i\times P_i})=\frac{\rho^D}{\mu(P_i)}$, and with left and right eigenfunctions $w_L:=\sqrt{\mu(P_i)}\cdot (v_L)_{\restriction P_i}$ and $w_R:=\sqrt{\mu(P_i)}\cdot(v_R)_{\restriction P_i}$. The rescaling by the factor $\sqrt{\mu(P_i)}$ is chosen so that we have $\langle w_L,w_R\rangle_{P_i}=1$ for the inner product $\langle \cdot,\cdot\rangle_{P_i}$ in the Hilbert space $L^2(P_i,\mu_{P_i})$.

By Proposition~\ref{prop:asymptoticAperiodic}, we have
\[
\left((\Gamma^D)_{\restriction P_i\times P_i}\right)^M(x,y)= \left(\frac{\rho^D}{\mu(P_i)}\right)^M w_R(x)w_L(y)+O(\beta^M)=\frac{\rho^{MD}}{\mu(P_i)^{M-1}}
v_R(x)v_L(y)+O(\beta^M)\;,
\]
for some $\beta\in(0,\frac{\rho^D}{\mu(P_i)})$. Thus, the statement follows by substituting into~\eqref{eq:SubstIntoMe}. Here, $\alpha_{i,i}=\sqrt[D]{\beta \mu(P_i)}<\rho$.

\underline{Next, we prove the case when $i\neq j$.} Suppose that $(x,y)\in P_i\times P_j$. Set $M^*:=M-1$ and $r:=D+j-i$. Define a function $F\in L^2(\Omega)$ by $F(z):=\Gamma^{M^*D}(x,z)$. This function is~0 if $z\not\in P_i$. If $z\in P_i$, then the previously established part applies. Thus, for the functions $G,H\in L^2(\Omega)$, 
\begin{equation*}
G(z):=
\begin{cases}
\rho^{M^* D} v_R(x) v_L(z)   &\mbox{if $z\in P_i$, and}\\
0                                           &\mbox{otherwise,}
\end{cases}
\end{equation*}
and $H:=F-G$,
we have $\|H\|_\infty= O(\alpha_{i,i}^{M^*D})$. Below, we express $\Gamma^{\ell}(x,y)$ using the operator $\Gamma^r$ acting on the function $F$,
\[
\Gamma^{\ell}(x,y)=\left((\Gamma^r)(F)\right)(y)=\left((\Gamma^r)\left(G+H\right)\right)(y)
=
\left((\Gamma^r)\left(G\right)\right)(y)
+
\left((\Gamma^r)\left(H\right)\right)(y)
\;,
\]
by linearity. We deal with both summands separately. The term $\left((\Gamma^r)\left(G\right)\right)(y)$ will lead to the main term of~\eqref{TE}, while the term $\left((\Gamma^r)\left(H\right)\right)(y)$ will be the error term.

Let us look at the term $\left((\Gamma^r)\left(G\right)\right)(y)$. The function $G$ is the $\rho^{M^* D} v_R(x)$-multiple of the eigenfunction $v_L$ on $P_i$ and is~0 on other cyclic sets. The periodicity of $\Gamma$ means, that when for some $f\in L^2(\Omega)$ we want to evaluate $\Gamma f$ at some cyclic set $P_k$, only the values of $f$ on the previous cyclic set $P_{k-1}$ are relevant. This demonstrates particularly nicely if $f$ coincides with some eigenfunction $f_\tau$ (whose eigenvalue is $\tau\in \C$) on $P_{k-1}$. In that case $(\Gamma f)_{\restriction P_k}=\tau f_{\restriction P_k}$. Of course, the relevant eigenfunction for us is $\rho^{M^* D} v_R(x)\ v_L$. Using the above iteratively $r$ many times with $k\equiv i+1,i+2,\ldots,i+r-1\mod D$, we have that $\left((\Gamma^r)\left(G\right)\right)$ is the $\rho^r$-multiple of $\rho^{M^* D} v_R(x)\ v_L$. That is, $\left((\Gamma^r)\left(G\right)\right)(y)=\rho^r\cdot\rho^{M^* D} v_R(x)\ v_L(y)=\rho^{\ell} v_R(x)\ v_L(y)$, as appears in~\eqref{TE}.

To deal with the term $\left((\Gamma^r)\left(H\right)\right)(y)$ is simple. We have $\| \Gamma\|_\infty\le 1$. That is, the operator $\Gamma$ is contractive with respect to the $L^\infty$-norm. Thus, $\|(\Gamma^r)\left(H\right)\|_\infty \le \| H\|_\infty = O(\alpha_{i,i}^{M^*D})=O(\alpha_{i,i}^{\ell})$.

\section*{Acknowledgements}
We thank Tomasz Kania and Vladimír Müller for their help regarding functional analysis.

\addcontentsline{toc}{section}{Index}
%INDEX \printindex

\bibliographystyle{acm}
\addcontentsline{toc}{section}{Bibliography}
\bibliography{references.bib}

\section*{Contact details}
\noindent\begin{tabular}{ll} 
\emph{Postal address:} & Institute of Computer Science of the Czech Academy of Sciences \\ 
 & Pod Vodárenskou věží 2 \\ 
 & 182~00, Prague \\ 
 & Czechia\\
& \\
\emph{Email:}&\texttt{hladky@cs.cas.cz}, \texttt{savicky@cs.cas.cz}
\end{tabular}
\appendix

\end{document}